\documentclass[12pt, a4]{article}
\usepackage{amsmath, amssymb, bm}
 \usepackage{amsthm}
 \usepackage{physics}
 \usepackage{float}
 \usepackage{multicol}
\newtheorem{conj}{Conjecture}

\newtheorem{Rem}{Remark}

\newcommand{\ba}{\begin{eqnarray}}
\newcommand{\ea}{\end{eqnarray}}
\newcommand{\ban}{\begin{eqnarray*}}
\newcommand{\ean}{\end{eqnarray*}}
\newcommand{\no}{\nonumber}

\def\d{{\partial}}
\newcommand{\mapright}[1]{%
\smash{\mathop{%
\hbox to 1.0cm{\rightarrowfill}}\limits^{#1}}}
\newcommand{\mapleft}[1]{%
\smash{\mathop{%
\hbox to 1.3cm{\leftarrowfill}}\limits^{#1}}}
\begin{document}
\title{
\begin{flushright}
  \begin{minipage}[b]{5em}
    \normalsize
    ${}$      \\
  \end{minipage}
\end{flushright}
{\bf Enumeration of Elliptic Curves via Elliptic Gromov-Witten Invariants of Four Dimensional Projective Fano Hypersurfaces}}
\author{Masao Jinzenji${}^{(1)}$, Ken Kuwata${}^{(2)}$ \\
\\
${}^{(1)}$\it Department of Mathematics, Okayama University \\
\it  Okayama, 700-8530, Japan\\
${}^{(2)}$\it Department of General Education \\
\it National Institute of Technology, Kagawa College \\
\it Chokushi, Takamatsu, 761-8058, Japan\\
\\
\it e-mail address: \it\hspace{0.3cm}${}^{(1)}$ pcj70e4e@okayama-u.ac.jp \\
\it\hspace{4.0cm}${}^{(2)}$ kuwata-k@t.kagawa-nct.ac.jp}

\maketitle
\begin{abstract}
In this paper, we propose a conjecture that clarifies the relationship between the number of degree $d$ elliptic curves in complex four-dimensional projective Fano hypersurfaces and their degree $d$ elliptic Gromov-Witten (GW) invariants. The elliptic GW invariants are computed using the elliptic virtual structure constants proposed in our previous works\footnote{Keywords; enumerative geometry; gromov-witten invariants; mirror symmetry; elliptic curve; projective Fano hypersurface.  2020 MATHEMATICS SUBJECT CLASSIFICATION; 14N10, 14N35.   }.
\end{abstract}

\section{Introduction}
\subsection{Our Motivation}
The genus $g$ and degree $d$ Gromov-Witten (GW) invariants of a complex $n$-dimensional K\"{a}hler manifold $X$ are generally expected to count the number of holomorphic curves of genus $g$ and degree $d$ in $X$ that satisfy the passing-through conditions imposed by the operator insertions. Therefore, these invariants are anticipated to be non-negative integers. However, there are many cases where they fail to be non-negative integers. As far as we know, the genus $0$ GW invariants of projective Fano manifolds are always non-negative integers. In the case of Calabi-Yau and general-type manifolds, however, genus $0$ GW invariants typically become rational numbers due to contributions from multiple cover maps. For Calabi-Yau manifolds, many studies analyze the contributions from multiple cover maps, often relating them to open Calabi-Yau manifolds obtained as vector bundles over $CP^{1}$ \cite{AM,Kont}. In contrast, for general-type manifolds, there are few works on contributions from multiple cover maps to genus $0$ GW invariants, with the exception of the preprint \cite{JNS} in 2004.

Now, let us turn to the case of genus $1$ GW invariants. The genus $1$ GW invariants of projective Fano surfaces (complex two-dimensional manifolds) are non-negative integers, at least for lower degrees. However, in the Calabi-Yau case, we know that the elliptic GW invariants include not only contributions from multiple cover maps from a genus $1$ curve but also degenerate contributions associated with genus $0$ GW invariants \cite{BCOV,P1}. Roughly speaking, the degenerate contributions arise from the boundary loci of the moduli space of stable maps from genus $1$ stable curves, where the stable curves degenerate into nodal curves containing a genus $0$ component. At these loci, the stable curves can be mapped into rational curves in the target manifold $X$, and the contributions are related to the genus $0$ GW invariants of $X$. It is well-known that degenerate contributions also appear in the case of complex three-dimensional Fano manifolds. Heuristically, this phenomenon is expected from the dimensional counting of the virtual dimension of the moduli space of stable maps. The virtual dimension of the moduli space of stable maps from stable curves to complex three-dimensional projective manifolds with fixed degree $d$ is independent of the genus of the stable curves. Hence, we can speculate that the genus $1$ GW invariants may contain degenerate contributions because the operator insertions required to balance the positive virtual dimension of the moduli space do not depend on the genus. A natural question thus arises: {\bf "Do degenerate contributions appear in the case of a projective Fano manifold whose complex dimension is greater than 3 ?"}

In this paper, motivated by the above question, we investigate the relationship between the elliptic GW invariants of four-dimensional projective Fano manifolds and the number of elliptic curves. In the case of complex four-dimensional K\"{a}hler manifolds, the virtual dimension of the moduli space of stable maps from genus $0$ stable curves is one complex dimension greater than that from genus $1$ stable curves. Therefore, it is not clear whether degeneration contributions appear, based on the naive speculation used in the preceding discussion. In the Calabi-Yau case, the work \cite{KP} by Klemm and Pandharipande shows that degenerate contributions do indeed appear, and they are related to genus $0$ GW invariants in a more complicated way than in the case of complex three-dimensional Calabi-Yau manifolds. We must note that their results rely on the computation of genus $1$ GW invariants using the BCOV conjecture \cite{BCOV}. On the other hand, in \cite{JK}, we proposed a conjecture that enables a kind of B-model computation (from the perspective of mirror symmetry) of genus $1$ GW invariants for projective hypersurfaces of arbitrary dimension. We therefore used our conjecture to evaluate genus $1$ GW invariants of four-dimensional projective Fano hypersurfaces. In the Fano case, the Virasoro conjecture \cite{virasoro1} and Getzler's equation \cite{Getzler} are also effective for this purpose. When applying these methods to evaluate genus $1$ invariants, one must use information from genus $0$ GW invariants with insertions of operators derived from primitive cohomology classes of the hypersurfaces. As was done in \cite{virasoro2}, these genus $0$ invariants were computed in the case of a cubic hypersurface in $CP^4$ using the associativity equation, since it has primitive classes only in the bidegree $(2,1)$ and $(1,2)$ sectors. However, in the four-dimensional case, a projective Fano hypersurface can have primitive classes in the bidegree $(3,1)$, $(2,2)$, and $(1,3)$ sectors. Therefore, we must overcome technical obstacles to compute genus $0$ invariants with these operator insertions using known techniques. Moreover, as we experienced in the computation of genus $1$ invariants of the cubic hypersurface using the Virasoro conjecture \cite{virasoro2}, the method requires outrageously complicated and lengthy computations. In contrast, our conjecture does not require any information regarding genus $0$ GW invariants with insertions of primitive classes and can be easily automated using computer software such as Maple and Mathematica.

After computing the genus $1$ GW invariants of the degree $k$ ($k=1,2,3,4,5$) hypersurface in $CP^5$ (denoted by $M_{6}^{k}$) up to degree $d=5$, we concluded that degeneration contributions do appear because the invariants often turn out to be negative rational numbers. We then investigated how these contributions are described by the genus $0$ GW invariants of the hypersurface. The guiding principles we used for this determination were:
\begin{enumerate}
    \item The number of elliptic curves of degree $d=1,2$ in $M_{6}^{k}$ that satisfy the passing-through conditions imposed by operator insertions is always $0$.
    \item The number of elliptic curves of degree $d=3$ in $M_{6}^{k}$ that satisfy the passing-through conditions imposed by operator insertions can be computed by applying the method presented in the preprint \cite{Katz} by S. Katz.
    \item The expected formula is considered an extension of the result presented in \cite{KP} to the Fano case.
\end{enumerate}
Using these hints, we constructed a conjectural formula that describes the degenerate contributions in terms of genus $0$ GW invariants, which will be presented in the next section. Heuristically, the genus $1$ GW invariants of degree $d$ are considered the sum of the degenerate contributions and the number of elliptic curves of the same degree satisfying the passing-through conditions. Therefore, if we subtract the degenerate contributions from the genus $1$ GW invariants, we expect to obtain non-negative integers. We denote the resulting number by $E_{d,a,b,c}$, and it turns out to be always a non-negative integer up to $d=5$. Specifically, $E_{d,a,b,c}=0$ if $d=1,2$, and it coincides with the number of elliptic curves computed by S. Katz's method if $d=3$.  With these results, we propose a conjectural method for counting elliptic curves in a four-dimensional projective Fano hypersurface via genus $1$ GW invariants.

\subsection{Definitions and Our Main Conjecture}

Let $M_N^k$ denote a degree $k$ hypersurface in $CP^{N-1}$. Let $h$ be the pull-back of the hyperplane class of $H^{1,1}(CP^{N-1},{\bf C})$ to $H^{1,1}(M_N^k,{\bf C})$ via the inclusion map $\iota: M_{N}^{k}\hookrightarrow CP^{N-1}$. We denote the genus $g$ and degree $d$ Gromov-Witten (GW) invariants of $M_N^k$ by $\langle\prod_{j=1}^{L}({\cal O}_{\alpha_{j}})^{m_{j}}\rangle_{g,d}$, where $\alpha_{j}$'s ($j=1,\cdots, L$) are a linear basis of $H^{*,*}(M_{N}^{k},{\bf C})$.
We omit the rigorous definition of the Gromov-Witten invariant here. In this paper, we only consider the following GW invariants:
\ba
\langle\prod_{a=0}^{N-2}({\cal O}_{h^a})^{m_{a}}\rangle_{0,d}\quad (g=0,1).
\label{type}
\ea
In the $g=0$ and $g=1$ cases, it is well-known that the GW invariants are non-zero only if the following conditions are satisfied:
\begin{align}
&\langle\prod_{a=0}^{N-2}({\cal O}_{h^a})^{m_{a}}\rangle_{0,d}\neq 0\;\Longrightarrow \sum_{a=0}^{N-2}m_{a}(a-1)=N-5+d(N-k),\;\;\\
&\langle\prod_{a=0}^{N-2}({\cal O}_{h^a})^{m_{a}}\rangle_{1,d}\neq 0\;\Longrightarrow \sum_{a=0}^{N-2}m_{a}(a-1)=d(N-k).\;\;
\label{sel}
\end{align}
For later use, we also introduce the {\bf K\"{a}hler (divisor) axiom} and the {\bf puncture axiom} of the GW invariants.
\ba
&& \langle{\cal O}_{h}\prod_{a=0}^{N-2}({\cal O}_{h^a})^{m_{a}}\rangle_{g,d}=d\langle\prod_{a=0}^{N-2}({\cal O}_{h^a})^{m_{a}}\rangle_{g,d}\quad (d\geq 1), \no\\
&& \langle{\cal O}_{h^{0}}\prod_{a=0}^{N-2}({\cal O}_{h^a})^{m_{a}}\rangle_{g,d}=0\hspace{3.7cm} (d\geq 1).
\label{ax}
\ea
If $d=0$, the genus $0$ GW invariants are non-vanishing only if the number of operator insertions $\sum_{a=0}^{N-2}m_{a}$ equals $3$, and the non-vanishing invariants are given by:
\ba
\langle{\cal O}_{h^a}{\cal O}_{h^b}{\cal O}_{h^{c}}\rangle_{0,0}=\int_{M_{N}^{k}}h^{a}\wedge h^{b}\wedge h^{c}=k\cdot\delta_{a+b+c,N-2}.
\ea
If $g=1$ and $d=0$, the only non-vanishing GW invariant is $\langle {\cal O}_{h}\rangle_{1,0}$, which is given by
\ba
\langle {\cal O}_{h}\rangle_{1,0}=-\frac{1}{24}\int_{M_{N}^{k}}h\wedge c_{N-3}(T^{\prime}M_{N}^{k}),
\ea
where $T^{\prime}M_{N}^{k}$ is the holomorphic tangent bundle of $M_{N}^{k}$ and $c_{N-3}(T^{\prime}M_{N}^{k})$ is its second-to-top Chern class.

From now on, we fix $N=6$. Hence, we only consider $M_{6}^{k}$, the four-dimensional hypersurface in $CP^{5}$. Considering the K\"{a}hler equation and the puncture equation, we can easily see that the non-trivial $g=0,1$ GW invariants with positive degrees are given by
\ba
\langle({\cal O}_{h^2})^{a}({\cal O}_{h^3})^{b}({\cal O}_{h^4})^{c}\rangle_{0,d},\quad 
\langle({\cal O}_{h^2})^{a}({\cal O}_{h^3})^{b}({\cal O}_{h^4})^{c}\rangle_{1,d} \quad (d\geq 1)
\label{type2}
\ea
From (\ref{sel}), all $g=1$ and $d\geq 1$ GW invariants (given by (\ref{type2})) vanish if $k>6$. In the $k=6$ case, the non-trivial genus $0$ and $1$ GW invariants are given by $\langle{\cal O}_{h^2}\rangle_{0,d}$ and $\langle*\rangle_{1,d}$ ($d\geq 1$), respectively ($*$ means no operator insertions). The relation between $\langle*\rangle_{1,d}$ and the number of elliptic curves of degree $d$ in $M_{6}^{6}$ was extensively studied in \cite{KP}. Therefore, we focus on the $k=1,2,3,4,5$ cases. In these cases, we evaluated the genus $0$ GW invariants in (\ref{type2}) using both the multi-point virtual constants \cite{JS} and the associativity equation \cite{KM}. We then evaluated the corresponding genus $1$ invariants using the method of elliptic virtual structure constants \cite{JK}\footnote{{\bf Text copies} of {\bf Mathematica programs} for these computations are available in {\bf ReserchGate homepage of Ken Kuwata}\cite{Kuwata}.}. Based on these numerical results and the three criteria introduced in Section 1, we propose the following conjecture:
\begin{conj}
Let $E_{d,a,b,c}$ be defined by the following relation:
\begin{align}
&\langle({\cal O}_{h^2})^{a}({\cal O}_{h^3})^{b}({\cal O}_{h^4})^{c}\rangle_{1,d}=-\frac{(k^2-6k+15)d-(6-k)}{24d}\langle({\cal O}_{h^2})^{a+1}({\cal O}_{h^3})^{b}({\cal O}_{h^4})^{c}\rangle_{0,d}\notag \\
&+\frac{1}{2}\sum_{\substack{0<d_1, d_2\\d_1+d_2=d}}\frac{2d-d_1d_2(6-k)}{24d}\sum_{l=0}^a\sum_{m=0}^b\sum_{n=0}^c\binom{a}{l}\binom{b}{m}\binom{c}{n}\notag \\
&\times\left(\frac{1}{k}\langle({\cal O}_{h^2})^{l}({\cal O}_{h^3})^{m}({\cal O}_{h^4})^{n}{\cal O}_{h^3}\rangle_{0,d_1}\langle {\cal O}_{h}({\cal O}_{h^2})^{a-l}({\cal O}_{h^3})^{b-m}({\cal O}_{h^4})^{c-n}\rangle_{0,d_2}\right.\notag \\&
\left.+\frac{1}{k}\langle({\cal O}_{h^2})^{l}({\cal O}_{h^3})^{m}({\cal O}_{h^4})^{n} {\cal O}_{h}\rangle_{0,d_1}\langle {\cal O}_{h^3}({\cal O}_{h^2})^{a-l}({\cal O}_{h^3})^{b-m}({\cal O}_{h^4})^{c-n}\rangle_{0,d_2}\right. \notag\\
&\left.+\frac{1}{k} \langle({\cal O}_{h^2})^{l}({\cal O}_{h^3})^{m}({\cal O}_{h^4})^{n} {\cal O}_{h^2}\rangle_{0,d_1}\langle {\cal O}_{h^2}({\cal O}_{h^2})^{a-l}({\cal O}_{h^3})^{b-m}({\cal O}_{h^4})^{c-n}\rangle_{0,d_2}\right)\notag\\
&+\frac{k-2}{48}a!b!c!\langle({\cal O}_{h^2})^{a}({\cal O}_{h^3})^{b}({\cal O}_{h^4})^{c}\rangle_{0,1}\delta_{k,5}\delta_{d,2}+E_{d,a,b,c}.
\label{main}
\end{align}
Here, $\binom{a}{l}$ denotes a binomial coefficient, and $\delta_{i,j}$ denotes the Kronecker delta.
Then, $E_{d,a,b,c}$ is expected to be an integer representing the number of elliptic curves.
\end{conj}
\begin{Rem}
The first and second Chern classes of $T^{\prime}M_{6}^{k}$ are given by
\ba
c_{1}(T^{\prime}M_{6}^{k})=(6-k)h,\;\;c_{2}(T^{\prime}M_{6}^{k})=(k^2-6k+15)h^2.
\ea
Note the appearance of the polynomials in $k$, $(6-k)$ and $(k^2-6k+15)$, in (\ref{main}).
\end{Rem}   

Heuristically, we expect that {\bf $E_{d,a,b,c}$ is the number of elliptic curves of degree $d$ in $M^{k}_{6}$ passing through $a$ $3$-planes, $b$ $2$-planes, and $c$ $1$-planes in $CP^5$}. In particular, for $d = 1, 2$, our numerical computation confirmed that $E_{d,a,b,c}$ always vanishes.
Furthermore, for $d = 3$, we confirmed that $E_{d,a,b,c}$ coincides with the number of elliptic curves satisfying the required passing-through conditions, as computed by applying the method presented in \cite{Katz} by S. Katz. We will explain the outline of S. Katz's method in Appendix A.

\section{Derivation of the Conjecture}
\subsection{The $d=1$ case}
We begin with the $d=1$ case. In this setting, it is known that no elliptic curves exist within the hypersurface. Consequently, the Gromov-Witten invariant $\langle({\cal O}_{h^{2}})^a({\cal O}_{h^{3}})^b({\cal O}_{h^{4}})^c\rangle_{1,1}$ consists solely of the degenerate contribution related to the genus $0$ and degree $1$ invariant $\langle({\cal O}_{h^{2}})^e({\cal O}_{h^{3}})^f({\cal O}_{h^{4}})^g\rangle_{0,1}$. 

In evaluating this degenerate contribution, the incidence conditions arising from operator insertions are effective. However, the virtual dimension of the genus $0$ moduli space exceeds that of the genus $1$ moduli space by $1$. We therefore speculate that the contribution can be expressed using the genus $0$ Gromov-Witten invariant $\langle({\cal O}_{h^{2}})^{a+1}({\cal O}_{h^{3}})^b({\cal O}_{h^{4}})^c\rangle_{0,1}$. Let us examine the data for these invariants for $M_{6}^{1}$, as shown in the tables below \footnote{For simplicity, we denote the genus $g$, degree $d$ Gromov-Witten invariants of $M_6^k$,\\ 
$\left< ( {\cal O}_{h^2})^{a} ( {\cal O}_{h^3})^{b} ( {\cal O}_{h^4})^{c}\right>_{g,d}$, by $N^g_{d,a,b,c}$.}.

\begin{table}[H]
\centering
\caption{$g=0, M_6^1, d=1$} 
\begin{tabular}{|c|l|l|}
\hline
d & (a,b,c) & $N^0_{d,a,b,c}$ \\
\hline
1 & (0,0,2) & $1$ \\ \hline
1 & (1,1,1) & $1$ \\ \hline
1 & (3,0,1) & $1$ \\ \hline
1 & (0,3,0) & $1$ \\ \hline
1 & (2,2,0) & $2$ \\ \hline
1 & (4,1,0) & $3$ \\ \hline
1 & (6,0,0) & $5$ \\ \hline
\end{tabular}
\end{table}

\begin{table}[H]
\centering
\caption{$g=1, M_6^1, d=1$} \label{g=1,N=6,k=1}
\begin{tabular}{|c|l|l|}
\hline
d & (a,b,c) & $N^1_{d,a,b,c}$ \\
\hline
1 & (0,1,1) & $-\frac{5}{24}$ \\ \hline
1 & (2,0,1) & $-\frac{5}{24}$ \\ \hline
1 & (1,2,0) & $-\frac{5}{12}$ \\ \hline
1 & (3,1,0) & $-\frac{5}{8}$ \\ \hline
1 & (5,0,0) & $-\frac{25}{24}$ \\ \hline
\end{tabular}
\end{table}

From this data, we readily identify the following relation for $k=1$:
\ba
\langle({\cal O}_{h^{2}})^a({\cal O}_{h^{2}})^b({\cal O}_{h^{2}})^c\rangle_{1,1}=-\frac{5}{24}\langle({\cal O}_{h^{2}})^{a+1}({\cal O}_{h^{3}})^b({\cal O}_{h^{4}})^c\rangle_{0,1}.
\label{k1d1}
\ea
Similar relations were found for other values of $k$:
\ba
&&\langle({\cal O}_{h^{2}})^a({\cal O}_{h^{2}})^b({\cal O}_{h^{2}})^c\rangle_{1,1}=-\frac{3}{24}\langle({\cal O}_{h^{2}})^{a+1}({\cal O}_{h^{3}})^b({\cal O}_{h^{4}})^c\rangle_{0,1}\quad (k=2), \no\\ 
&&\langle({\cal O}_{h^{2}})^a({\cal O}_{h^{2}})^b({\cal O}_{h^{2}})^c\rangle_{1,1}=-\frac{3}{24}\langle({\cal O}_{h^{2}})^{a+1}({\cal O}_{h^{3}})^b({\cal O}_{h^{4}})^c\rangle_{0,1}\quad (k=3),\no\\
&&\langle({\cal O}_{h^{2}})^a({\cal O}_{h^{2}})^b({\cal O}_{h^{2}})^c\rangle_{1,1}=-\frac{5}{24}\langle({\cal O}_{h^{2}})^{a+1}({\cal O}_{h^{3}})^b({\cal O}_{h^{4}})^c\rangle_{0,1}\quad (k=4),\no\\
&&\langle({\cal O}_{h^{2}})^a({\cal O}_{h^{2}})^b({\cal O}_{h^{2}})^c\rangle_{1,1}=-\frac{9}{24}\langle({\cal O}_{h^{2}})^{a+1}({\cal O}_{h^{3}})^b({\cal O}_{h^{4}})^c\rangle_{0,1}\quad (k=5).
\label{okd1}
\ea
To unify these results, we assumed the following ansatz: 
\ban
\langle({\cal O}_{h^{2}})^a({\cal O}_{h^{2}})^b({\cal O}_{h^{2}})^c\rangle_{1,1}=-\frac{ak^2+bk+c}{24}\langle({\cal O}_{h^{2}})^{a+1}({\cal O}_{h^{3}})^b({\cal O}_{h^{4}})^c\rangle_{0,1}.
\ean
Solving for the coefficients, we obtained the general relation for $M_{6}^{k}$:
\ba
\langle({\cal O}_{h^{2}})^a({\cal O}_{h^{2}})^b({\cal O}_{h^{2}})^c\rangle_{1,1}=-\frac{k^2-5k+9}{24}\langle({\cal O}_{h^{2}})^{a+1}({\cal O}_{h^{3}})^b({\cal O}_{h^{4}})^c\rangle_{0,1}.
\ea

\subsection{Partial Results for the $d=2$ and $d=3$ cases}

Next, we considered the $d=2$ case. In this case, there are no elliptic curves in $M_{6}^{k}$. Therefore, we initially assumed a relation of the form:
\ba
\langle({\cal O}_{h^{2}})^a({\cal O}_{h^{2}})^b({\cal O}_{h^{2}})^c\rangle_{1,2} = -\frac{\alpha(k,2)}{24}\langle({\cal O}_{h^{2}})^{a+1}({\cal O}_{h^{3}})^b({\cal O}_{h^{4}})^c\rangle_{0,2}.
\ea 
However, for cases other than $k=2$, additional contributions appear to exist, and the above ansatz does not hold. For the specific case $k=2$, we identified the following relation:
\ba
\langle({\cal O}_{h^{2}})^a({\cal O}_{h^{2}})^b({\cal O}_{h^{2}})^c\rangle_{1,2} = -\frac{5}{24}\langle({\cal O}_{h^{2}})^{a+1}({\cal O}_{h^{3}})^b({\cal O}_{h^{4}})^c\rangle_{0,2} \quad (k=2).
\ea 

We then turned to the $d=3$ case. Here, there is a non-trivial contribution to $\langle({\cal O}_{h^{2}})^a({\cal O}_{h^{2}})^b({\cal O}_{h^{2}})^c\rangle_{1,3}$ from the number of elliptic curves, denoted by $E_{3,a,b,c}$. 
At the early stages of this study, a general method for evaluating $E_{3,a,b,c}$ for arbitrary $M_{6}^{k}$ had not yet been established. However, for $k=3$, it is possible to count $E_{3,a,b,c}$ directly because a degree $3$ elliptic curve in $M_{6}^{3}$ is always given by the intersection of a $2$-plane in $CP^{5}$ and $M_{6}^{3}$. 

Let $H_{a,b,c}$ be the number of $2$-planes that pass through $a$ $3$-planes, $b$ $2$-planes, and $c$ $1$-planes in $CP^5$. This can be readily evaluated using Schubert calculus on the Grassmannian $Gr(3,6)$. Then, $E_{3,a,b,c}$ for $M_{6}^{3}$ is given by:
\ba
E_{3,a,b,c} = 3^{a+b+c} H_{a,b,c} \quad (k=3). 
\ea
Using the data obtained through this method, we found that the following relation holds for the $k=3, d=3$ case:
\ba
\langle({\cal O}_{h^{2}})^a({\cal O}_{h^{2}})^b({\cal O}_{h^{2}})^c\rangle_{1,3} = -\frac{5}{24}\langle({\cal O}_{h^{2}})^{a+1}({\cal O}_{h^{3}})^b({\cal O}_{h^{4}})^c\rangle_{0,3} + E_{3,a,b,c} \quad (k=3).
\ea

\subsection{Incorporating the Results of Klemm and Pandharipande}  
In the $d=2$ and $d=3$ cases, we had not yet determined the additional contributions expected to appear, except for the specific cases mentioned in the previous subsection. To address this, we consulted the work of Klemm and Pandharipande \cite{KP}, which provides the solution for the four-dimensional Calabi-Yau hypersurface $M_{6}^{6}$. 

In the $M_{6}^{6}$ case, the non-trivial genus $0$ and genus $1$ Gromov-Witten invariants of degree $d$ are:
 \ban
 \langle{\cal O}_{h^2}\rangle_{0,d}, \quad \langle * \rangle_{1,d},
\ean
where the symbol $*$ denotes the absence of operator insertions. Since $M_{6}^{6}$ is a Calabi-Yau manifold with vanishing first Chern class, these are the only non-vanishing invariants for these genera. To isolate the contributions of rational curves from multiple cover maps, we define the enumerative invariants $n_{0,d}$--representing the number of rational curves of degree $d$ passing through a $3$-plane in $CP^{5}$--via the following relation:
\ba
\sum_{d=1}^\infty \langle{\cal O}_{h^2}\rangle_{0,d}q^d = \sum_{d=1}^\infty n_{0,d} \sum_{s=1}^\infty \frac{q^{sd}}{s^2}. \label{nrc}
\ea
The result in \cite{KP} for the Calabi-Yau hypersurface $M_{6}^{6}$ is given by:
\ba
\sum_{d=1}^\infty \langle 1 \rangle_{1,d}q^d &=& \sum_{d=1}^\infty n_{1,d} \sum_{s=1}^\infty \frac{\sigma(s)q^{sd}}{s} + \frac{1}{24} \sum_{d=1}^\infty c_2(M^6_6) n_{0,d} \log(1-q^d) \no \\
&& - \frac{1}{24} \sum_{d_1,d_2=1}^\infty m_{d_1,d_2} \log(1-q^{d_1+d_2}), \label{nec}
\ea
where $n_{1,d}$ is the number of elliptic curves of degree $d$ in $M_{6}^{6}$ (corresponding to $E_{d,a,b,c}$ in our notation), $c_2(M^6_6) = 15$ is the second Chern class coefficient, and $m_{d_1,d_2}$ is a rational number defined by the following recursive rules:

\begin{enumerate}
    \item $m_{d_1,d_2} = m_{d_2,d_1}$
    \item $m_{d_1,d_2} = 0$ if $d_1 \leq 0$ or $d_2 \leq 0$.
    \item If $d_1 \neq d_2$:
    \ban
    m_{d_1,d_2} = \frac{n_{0,d_1}n_{0,d_2}}{k} + m_{d_1, d_2-d_1} + m_{d_1-d_2, d_2}
    \ean
    \item If $d_1 = d_2 = d$:
    \ban
    m_{d,d} = c_2(M_6^6)n_{0,d} + \frac{n_{0,d}n_{0,d}}{k} - \sum_{d_1+d_2=d} m_{d_1,d_2}
    \ean
\end{enumerate}

Formula (\ref{nec}) is complex due to multiple cover contributions arising from the vanishing first Chern class. However, in the Fano case, these contributions can often be ignored or simplified. We therefore hypothesized that a simpler structure exists and constructed the following ansatz for the Fano hypersurfaces:
\ba
&&\langle({\cal O}_{h^2})^{a}({\cal O}_{h^3})^{b}({\cal O}_{h^4})^{c}\rangle_{1,d} = E_{d,a,b,c} - \frac{\alpha(k,d)}{24}\langle({\cal O}_{h^2})^{a+1}({\cal O}_{h^3})^{b}({\cal O}_{h^4})^{c}\rangle_{0,d} \no \\
&&+ \frac{\beta(k,d)}{24} \sum_{\substack{d_1+d_2=d \\ d_1, d_2 > 0}} \sum_{\alpha=1}^{3} \sum_{l=0}^a \sum_{m=0}^b \sum_{n=0}^c \binom{a}{l}\binom{b}{m}\binom{c}{n} \no \\
&&\times \frac{1}{k} \langle({\cal O}_{h^2})^{l}({\cal O}_{h^3})^{m}({\cal O}_{h^4})^{n}{\cal O}_{h^{\alpha}}\rangle_{0,d_1} \langle {\cal O}_{h^{4-\alpha}}({\cal O}_{h^2})^{a-l}({\cal O}_{h^3})^{b-m}({\cal O}_{h^4})^{c-n}\rangle_{0,d_2}. 
\label{anz2}
\ea

The three terms on the right-hand side of (\ref{anz2}) correspond to the three components in (\ref{nec}). The second term, involving $\alpha(k,d)$, was partially determined from the $d=1,2$ cases as:
\ba
\alpha(k,1) = k^2-5k+9 \quad (1 \leq k \leq 6), \quad \alpha(2,2)=5, \quad \alpha(3,3)=5.
\ea 
For $M_6^6$, (\ref{nec}) implies $\alpha(6,1) = c_2(M_6^6)$. Given that the first and second Chern classes of $M_6^k$ are $(6-k)h$ and $(k^2-6k+15)h^2$ respectively, we identified $\alpha(k,1)$ as the combination $(k^2-6k+15) - (6-k)$. Following the logic that $\alpha(6,d) = 15$ implied by (\ref{nec}), we generalized the ansatz to:
\ban
\alpha(k,d) = (k^2-6k+15) - \gamma(d)(6-k).
\ean
Using $\alpha(2,2)=5$ and $\alpha(3,3)=5$, we found $\gamma(2)=1/2$ and $\gamma(3)=1/3$, leading to our final form:
\ba
\alpha(k,d) = (k^2-6k+15) - \frac{6-k}{d}.
\label{anz-alpha}
\ea 

The third term in (\ref{anz2}) is motivated by the boundary of the moduli space corresponding to nodal curves. Applying the splitting axiom, we assumed the singularity corresponds to the sum over operators $\sum \frac{1}{k}|{\cal O}_{h^{\alpha}}\rangle \langle{\cal O}_{h^{4-\alpha}}|$. Numerically testing this for $d=2$ (where $E_{2,a,b,c}=0$), we found:
\ban
\beta(1,2)=-1/4, \quad \beta(2,2)=0, \quad \beta(3,2)=1/4, \quad \beta(4,2)=1/2,
\ean
which suggests $\beta(k,2) = (k-2)/4$. For $k=5$, an exceptional extra term:
\ba
\frac{1}{16}a!b!c!\langle({\cal O}_{h^2})^{a}({\cal O}_{h^3})^{b}({\cal O}_{h^4})^{c}\rangle_{0,1},
\ea       
appeared, likely because the virtual dimensions of the genus $1$ degree $2$ and genus $0$ degree $1$ moduli spaces coincide.

For $d \geq 3$, we required that all $E_{d,a,b,c}$ derived from the ansatz be non-negative integers. This led to the refined ansatz:
\ba
&&\langle({\cal O}_{h^2})^{a}({\cal O}_{h^3})^{b}({\cal O}_{h^4})^{c}\rangle_{1,d}=E_{d,a,b,c}-\frac{\alpha(k,d)}{24}\langle({\cal O}_{h^2})^{a+1}({\cal O}_{h^3})^{b}({\cal O}_{h^4})^{c}\rangle_{0,d}\no \\
&&+\frac{1}{48d}\sum_{\substack{0<d_1, d_2\\d_1+d_2=d}}\beta(k,(d_{1},d_{2}))\sum_{\alpha=1}^{3}\sum_{l=0}^a\sum_{m=0}^b\sum_{n=0}^c\binom{a}{l}\binom{b}{m}\binom{c}{n}\no \\
&&\times\frac{1}{k}\langle({\cal O}_{h^2})^{l}({\cal O}_{h^3})^{m}({\cal O}_{h^4})^{n}{\cal O}_{h^{\alpha}}\rangle_{0,d_1}\langle {\cal O}_{h^{4-\alpha}}({\cal O}_{h^2})^{a-l}({\cal O}_{h^3})^{b-m}({\cal O}_{h^4})^{c-n}\rangle_{0,d_2}\quad (d\geq 3),\no\\ 
\label{anz4}
\ea

Our numerical experiments yielded the following values for $\beta(k,(d_1,d_2))$:

\begin{table}[H]
\centering
\caption{Numerical values for $\beta(k,(d_1,d_2))$} \label{beta}
\begin{tabular}{|c|l||c|l|}
\hline
$(d_1,d_2)$ & $\beta(k,(d_1,d_2))$ & $(d_1,d_2)$ & $\beta(k,(d_1,d_2))$ \\
\hline
(1,1) & $k-2$   & (1,4) & $4k-14$ \\ \hline
(1,2) & $2k-6$  & (2,2) & $4k-16$ \\ \hline
(1,3) & $3k-10$ & (2,3) & $6k-26$ \\ \hline
\end{tabular}
\end{table}

This allowed us to generalize the coefficient as $\beta(k,(d_1,d_2)) = 2d - d_1 d_2(6-k)$, leading to the final conjecture. In the tables below, we present our numerical data and the resulting $E_{d,a,b,c}$'s in the case of $M_{6}^{5}$. Subsequent validation using the method of S. Katz (Appendix A) confirmed the $d=3$ case.

\begin{table}[H]
\centering
\caption{$g=0,M_6^5$} \label{g=0,N=6,k=5}
\begin{tabular}{|c|l|l|l|l|l|}
\hline
d&(a,b,c)& $N^0_{d,a,b,c}$  \\
\hline
1&(0,1,0) & $3250$\\ \hline
1&(2,0,0) & $6125$\\ \hline
2&(0,0,1) & $247500$\\ \hline
2&(1,1,0)&$3718750$ \\ \hline
2&(3,0,0)&$17406875$ \\ \hline
3&(1,0,1) & $659250000$\\ \hline
3&(0,2,0) & $2700512500$\\ \hline
3&(2,1,0) & $19190225000$\\ \hline
3&(4,0,0) & $150549428125$\\ \hline
4&(0,1,1) & $529823250000$\\ \hline
4&(2,0,1) & $4729124250000$\\ \hline
4&(1,2,0) & $20567866625000$\\ \hline
4&(3,1,0) & $209913851312500$\\ \hline
4&(5,0,0) & $2337181124531250$\\ \hline
5&(0,0,2) & $109236016800000$\\ \hline
5&(1,1,1) & $5092187634000000$\\ \hline
5&(3,0,1) & $61252356251250000$\\ \hline
5&(0,3,0) & $21811124012125000$\\ \hline
5&(2,2,0) & $274577525136875000$\\ \hline
5&(4,1,0) & $3688919538904687500$\\ \hline
5&(6,0,0) & $53412041211701171875$\\ \hline
\end{tabular}
\end{table}
\begin{table}[H]
\centering
\caption{$g=1,M_6^5$} \label{g=1,N=6,k=5}
\begin{tabular}{|c|l|l|l|l|l|}
\hline
d&(a,b,c)& $N^1_{d,a,b,c}$&$E_{d,a,b,c}$  \\
\hline
1&(1,0,0) & $-\frac{18375}{8}$&$0$\\ \hline
2&(0,1,0) & $-\frac{8038625}{6}$&$0$\\ \hline
2&(2,0,0)&$-\frac{148120375}{24}$ &$0$\\ \hline
3&(0,0,1) & $-246708750$&$947500$\\ \hline
3&(1,1,0) & $-7059310625$&$14139375$\\ \hline
3&(3,0,0) & $-\frac{1319375639375 }{24} $&$93667500$\\ \hline
4&(1,0,1) & $-1762031831250$&$34818450000$\\ \hline
4&(0,2,0) & $-\frac{23015099134375}{3}$&$111635590625$\\ \hline
4&(2,1,0) & $-\frac{155974421459375}{2}$&$1070162215625$\\ \hline
4&(4,0,0) & $-\frac{5194583816496875}{6}$&$10767522628125$\\ \hline
5&(0,1,1) & $-1870863418500000$&$97094656425000$\\ \hline
5&(2,0,1) & $-22353636563531250 $&$1193292296250000$\\ \hline
5&(1,2,0) & $-\frac{303154762779953125}{3}$&$4376885918562500$\\ \hline
5&(3,1,0) & $-\frac{8130487814000265625}{6}$&$56589153313015625$\\ \hline
5&(5,0,0) & $-\frac{470853639798192671875}{24}$&$767616850277828125$\\ \hline
\end{tabular}
\end{table}

\newpage

\appendix

\section{Enumeration of Degree $3$ Elliptic Curves in $M_{6}^{k}$ ($k=1,2,3,4,5$)}
It is well-known that any degree $3$ elliptic curve in $CP^{N-1}$ is given as a **plane cubic**. Therefore, we first consider $Gr(3,N)$, the moduli space of $2$-planes in $CP^{N-1}$. A cubic curve in a plane (biholomorphic to $CP^2$) is determined by its defining equation, a degree $3$ homogeneous polynomial in three homogeneous coordinates of $CP^2$. Thus, the moduli space of degree $3$ elliptic curves in $CP^{N-1}$ is given by the fiber space $\pi_{P}:P(S^{3}U^{*})\rightarrow Gr(3,N)$, where $\pi_U:U\rightarrow Gr(3,N)$ is the tautological rank $3$ vector bundle over $Gr(3,N)$. 

Next, we construct the moduli space of elliptic curves in $M_{6}^{k}$ for $k$ varying from $1$ to $5$.
\begin{itemize}
    \item {\bf Case $k=1$}: Since $M_{6}^{1}$ is biholomorphic to $CP^4$, the moduli space is $\pi_{P}:P(S^{3}U^{*})\rightarrow Gr(3,5)$.
    \item {\bf Case $k=2$}: For a degree $3$ elliptic curve in $CP^5$ to be contained in $M_{6}^{2}$, the $2$-plane that contains the elliptic curve must itself be contained in $M_{6}^{2}$. Hence, the moduli space is given by $\pi_{P}^{-1}(PD(c_{\text{top}}(S^{2}U^{*})))$, where $PD(\alpha)$ is the submanifold of $Gr(3,6)$ that is Poincar\'{e} dual to the cohomology class $\alpha\in H^{*}(Gr(3,6))$, and $\pi_{P}: P(S^{3}U^{*})\rightarrow Gr(3,6)$ is the projection of the fiber space.
    \item {\bf Cases $k=3, 4, 5$}: For a degree $3$ elliptic curve in $CP^5$ to be contained in $M_{6}^{k}$, the restriction of the defining equation of $M_{6}^{k}$ (a homogeneous degree $k$ polynomial in $6$ homogeneous coordinates of $CP^5$) to the $2$-plane containing the elliptic curve must be divisible by the defining equation of the elliptic curve. Let $\pi_{S}:S_{P}\rightarrow P(S^{3}U^{*})$ be the tautological line bundle of the projectivization $P(S^{3}U^{*})$. Using this line bundle and the aforementioned condition, the moduli space of degree $3$ elliptic curves in $M_{6}^{k}$ is given by $PD(c_{\text{top}}(S^{k}U^{*}/(S^{k-3}U^{*}\otimes S_{P})))$.
\end{itemize}

Next, we briefly discuss how to express the condition that the elliptic curve passes through $PD(h^a)\subset CP^{N-1}$ in terms of $H^{*}(P(S^{3}U^{*}))$. For this purpose, we introduce the moduli space $\mathcal{M}$ of degree $3$ elliptic curves with one marked point. Let $\pi_{F}:\mathcal{M}\rightarrow P(S^{3}U^{*})$ be the forgetful map that omits the marked point, and let $ev: \mathcal{M}\rightarrow CP^{N-1}$ be the evaluation map. With these setups, the passing-through condition is represented by $\pi_{F*}(ev^{*}(h^a))\in H^{a-1}(P(S^{3}U^{*}))$, where $\pi_{F*}:H^{m}(\mathcal{M})\rightarrow H^{m-1}(P(S^{3}U^{*}))$ denotes fiber integration (or push-forward) by $\pi_{F}$. For details on the construction of $\mathcal{M}$, we recommend readers refer to \cite{Katz}. The cohomology class $\pi_{F*}(ev^{*}(h^a))$ is determined using the well-known projection formula \cite{bott-tu}, but we omit the computational details as they are highly technical. Instead, we present the explicit form of $\pi_{F*}(ev^{*}(h^a))$ at the end of this appendix.

Finally, we introduce the generators and relations of $H^{*}(P(S^{3}U^{*}))$ for explicit computation. Let $c_{i}$ ($i=1,2,3$) be the $i$-th Chern class of the vector bundle $U^{*}$. $H^{*}(Gr(3,N))$ is generated by these classes, and the relations are given as follows. Let $h_{i}(c_{1},c_{2},c_{3})$ be the weighted homogeneous polynomial of degree $i$ defined by
\ba
\frac{1}{1-c_{1}t+c_{2}t^2-c_{3}t^3}=1+\sum_{i=1}^{\infty}h_{i}t^{i}.
\ea
Then, we have
\ba
H^{*}(Gr(3,N))\simeq
{\bf C}[c_{1},c_{2},c_{3}]/(h_{N-2},h_{N-1},h_{N}),
\ea 
where $(h_{N-2},h_{N-1},h_{N})$ is the ideal of ${\bf C}[c_{1},c_{2},c_{3}]$ generated by $h_{N-2}, h_{N-1}$, and $h_{N}$. We now describe the ring structure of $H^{*}(P(S^{3}U^{*}))$. For this purpose, we introduce another generator $z$, which plays the role of the hyperplane class of the fiber projective space. Since $S^3U^{*}$ is a rank $10$ vector bundle over $Gr(3,N)$, $z$ satisfies the relation
\ba
R(z,c_{1},c_{2},c_{3}):=z^{10}+\sum_{j=1}^{10}c_{i}(S^{3}U^{*})z^{10-i}=0.
\ea
Hence, the ring structure of $H^{*}(P(S^{3}U^{*}))$ is given by
\ba
H^{*}(P(S^{3}U^{*}))\simeq {\bf C}[c_{1},c_{2},c_{3},z]/(h_{N-2},h_{N-1},h_{N},R).
\ea
The integration rule for cohomology elements of $H^{*}(P(S^{3}U^{*}))$ is fixed by the following equality:
\ba
\int_{P(S^{3}U^{*})}z^{9}\cdot(c_{3})^{N-3}=1.
\ea
The explicit forms of $\pi_{F*}(ev^{*}(h^a))$ for $a=2,3,4$ are given as follows:
\ban
\pi_{F*}(ev^{*}(h^2))&=&3c_{1}+z=:G_{1}\\
\pi_{F*}(ev^{*}(h^3))&=&3c_{1}^2-3c_{2}+c_{1}z=:G_{2}\\
\pi_{F*}(ev^{*}(h^4))&=&3c_{1}^3-6c_{1}c_{2}+3c_{3}+(c_{1}^2-c_{2})z=:G_{3}.
\ean 
These expressions do not depend on $N$.
Let $\tilde{E}_{3,a,b,c}$ be the number of degree $3$ elliptic curves in $M_{6}^{k}$ that corresponds to $E_{3,a,b,c}$ in our conjecture. 
With these preparations, $\tilde{E}_{3,a,b,c}$ is explicitly computed by the following formulae:
\ban
\tilde{E}_{3,a,b,c}&=&\int_{P(S^{3}U^{*})}G_{1}^aG_{2}^bG_{3}^c  \quad (k=1,\;N=5),\\
\tilde{E}_{3,a,b,c}&=&\int_{P(S^{3}U^{*})}c_{\text{top}}(S^{2}U^{*})G_{1}^aG_{2}^bG_{3}^c  \quad (k=2,\;N=6),\\
\tilde{E}_{3,a,b,c}&=&\int_{P(S^{3}U^{*})}c_{\text{top}}(S^{k}U^{*}/(S^{k-3}U^{*}\otimes S_{P}))G_{1}^aG_{2}^bG_{3}^c  \quad (k=3,4,5,\;N=6).
\ean  
We find that $\tilde{E}_{3,a,b,c}$ indeed coincides with $E_{3,a,b,c}$
\footnote{A \textbf{Text copy} of the \textbf{Maple worksheet} used for computing $\tilde{E}_{3,a,b,c}$ of $M_{6}^{5}$ is available on \textbf{M. Jinzenji's ResearchGate homepage} \cite{JR}.}.
In the $k=5$ case, we obtain the following results:
\ba
\tilde{E}_{3,0,0,1}=947500,\quad \tilde{E}_{3,1,1,0}=14139375,\quad \tilde{E}_{3,3,0,0}=93667500.
\ea

\newpage
\section{Lists of $g = 0, 1$ Gromov-Witten Invariants and 
$E_{d,a,b,c}$ for $M^k_6$ ($k = 1, 2, 3, 4$) up to $d = 5$}
\begin{multicols}{2}
\begin{table}[H]
\centering
\caption{$g=0,M_6^1,d=1,2$} \label{g=0,N=6,k=2}
\begin{tabular}{|c|l|l|l|l|l|}
\hline
d&(a,b,c)& $N^0_{d,a,b,c}$  \\
\hline
1&(0,0,2) & $1$\\ \hline
1&(1,1,1) & $1$\\ \hline
1&(3,0,1) & $1$\\ \hline
1&(0,3,0) & $1$\\ \hline
1&(2,2,0) & $2$\\ \hline
1&(4,1,0) & $3$\\ \hline
1&(6,0,0) & $5$\\ \hline
2&(0,1,3) & $0$\\ \hline
2&(2,0,3) & $1$\\ \hline
2&(1,2,2) & $1$\\ \hline
2&(3,1,2) & $4$\\ \hline
2&(5,0,2) & $11$\\ \hline
2&(0,4,1) & $2$\\ \hline
2&(2,3,1) & $6$\\ \hline
2&(4,2,1) & $21$\\ \hline
2&(6,1,1) & $67$\\ \hline
2&(8,0,1) & $219$\\ \hline
2&(1,5,0) & $10$\\ \hline
2&(3,4,0) & $36$\\ \hline
2&(5,3,0) & $132$\\ \hline
2&(7,2,0) & $473$\\ \hline
2&(9,1,0) & $1734$\\ \hline
2&(11,0,0) & $6620$\\ \hline
\end{tabular}
\end{table}

\begin{table}[H]
\centering
\caption{$g=0,M_6^1,d=3$} \label{g=0,N=6,k=1,d=3}
\begin{tabular}{|c|l|l|l|l|l|}
\hline
d&(a,b,c)& $N^0_{d,a,b,c}$  \\
\hline
3&(1,0,5) & $0 $\\ \hline
3&(0,2,4) & $1$\\ \hline
3&(2,1,4) & $5$\\ \hline
3&(4,0,4) & $30$\\ \hline
3&(1,3,3) & $9$\\ \hline
3&(3,2,3) & $45$\\ \hline
3&(5,1,3) & $225$\\ \hline
3&(7,0,3) & $1011 $\\ \hline
3&(0,5,2) & $16$\\ \hline
3&(2,4,2) & $76$\\ \hline
3&(4,3,2) & $385$\\ \hline
3&(6,2,2) & $1931$\\ \hline
3&(8,1,2) & $9386$\\ \hline
3&(10,0,2) & $45954$\\ \hline
3&(1,6,1) & $128 $\\ \hline
3&(3,5,1) & $664$\\ \hline
3&(5,4,1) & $3512 $\\ \hline
3&(7,3,1) & $18469$\\ \hline
3&(9,2,1) & $96548$\\ \hline
3&(11,1,1) & $511012 $\\ \hline
3&(13,0,1) & $2770596$\\ \hline
3&(0,8,0) & $188 $\\ \hline
3&(2,7,0) & $1108$\\ \hline
3&(4,6,0) & $6216$\\ \hline
3&(6,5,0) & $34780$\\ \hline
3&(8,4,0) & $194024$\\ \hline
3&(10,3,0) & $1085892$\\ \hline
3&(12,2,0) & $6165822$\\ \hline
3&(14,1,0) & $35806494$\\ \hline
3&(16,0,0) & $213709980$\\ \hline
\end{tabular}
\end{table}

\begin{table}[H]
\centering
\caption{$g=0,M_6^1,d=4$ part 1} \label{g=0,N=6,k=1,d=4,1}
\begin{tabular}{|c|l|l|l|l|l|}
\hline
d&(a,b,c)& $N^0_{d,a,b,c}$  \\
\hline
4&(0,0,7) & $1 $\\ \hline
4&(1,1,6) & $9 $\\ \hline
4&(3,0,6) & $61$\\ \hline
4&(0,3,5) & $14$\\ \hline
4&(2,2,5) & $107$\\ \hline
4&(4,1,5) & $732$\\ \hline
4&(6,0,5) & $4830$\\ \hline
4&(1,4,4) & $178$\\ \hline
4&(3,3,4) & $1218$\\ \hline
4&(5,2,4) & $8133$\\ \hline
4&(7,1,4) & $52507$\\ \hline
4&(9,0,4) & $324764$\\ \hline
4&(0,6,3) & $320$\\ \hline
4&(2,5,3) & $2056$\\ \hline
4&(4,4,3) & $13962$\\ \hline
4&(6,3,3) & $94104$\\ \hline
4&(8,2,3) & $622980$\\ \hline
4&(10,1,3) & $4063860$\\ \hline
4&(12,0,3) & $26578256$\\ \hline
4&(1,7,2) & $3516$\\ \hline
4&(3,6,2) & $23968$\\ \hline
4&(5,5,2) & $166936$\\ \hline
4&(7,4,2) & $1159218$\\ \hline
4&(9,3,2) & $7990720$\\ \hline
4&(11,2,2) & $54948346$\\ \hline
4&(13,1,2) & $380720598$\\ \hline
4&(15,0,2) & $2679044142$\\ \hline
\end{tabular}
\end{table}
\begin{table}[H]
\centering
\caption{$g=0,M_6^1,d=4$ part 2} \label{g=0,N=6,k=1,d=4,2}
\begin{tabular}{|c|l|l|l|l|l|}
\hline
d&(a,b,c)& $N^0_{d,a,b,c}$  \\
\hline
4&(0,9,1) & $5552$\\ \hline
4&(2,8,1) & $40492$\\ \hline
4&(4,7,1) & $291632$\\ \hline
4&(6,6,1) & $2110864$\\ \hline
4&(8,5,1) & $15251816$\\ \hline
4&(10,4,1) & $110031632$\\ \hline
4&(12,3,1) & $796460052$\\ \hline
4&(14,2,1) & $5823161346$\\ \hline
4&(16,1,1) & $43242657488$\\ \hline
4&(18,0,1) & $327439797532$\\ \hline
4&(1,10,0) & $63740$\\ \hline
4&(3,9,0) & $493976$\\ \hline
4&(5,8,0) & $3748804$\\ \hline
4&(7,7,0) & $28346212$\\ \hline
4&(9,6,0) & $213984472$\\ \hline
4&(11,5,0) & $1617593360$\\ \hline
4&(13,4,0) & $12302188692$\\ \hline
4&(15,3,0) & $94605276228$\\ \hline
4&(17,2,0) & $738764469204$\\ \hline
4&(19,1,0) & $5876564125104$\\ \hline
4&(21,0,0) & $47723447905060$\\ \hline
\end{tabular}
\end{table}
\begin{table}[H]
\centering
\caption{$g=0,M_6^1,d=5$ part 1} \label{g=0,N=6,k=1,d=5,1}
\begin{tabular}{|c|l|l|l|l|l|}
\hline
d&(a,b,c)& $N^0_{d,a,b,c}$  \\
\hline
5&(0,1,8) & $10$\\ \hline
5&(2,0,8) & $161 $\\ \hline
5&(1,2,7) & $246$\\ \hline
5&(3,1,7) & $2390$\\ \hline
5&(5,0,7) & $20670$\\ \hline
5&(0,4,6) & $432$\\ \hline
5&(2,3,6) & $3915$\\ \hline
5&(4,2,6) & $34180$\\ \hline
5&(6,1,6) & $284685$\\ \hline
5&(8,0,6) & $2269330$\\ \hline
5&(1,5,5) & $6700$\\ \hline
5&(3,4,5) & $57200$\\ \hline
5&(5,3,5) & $484345$\\ \hline
5&(7,2,5) & $4001415$\\ \hline
5&(9,1,5) & $32175350$\\ \hline
5&(11,0,5) & $252923350$\\ \hline
5&(0,7,4) & $11980$\\ \hline
5&(2,6,4) & $97660$\\ \hline
5&(4,5,4) & $829884$\\ \hline
5&(6,4,4) & $7042024$\\ \hline
5&(8,3,4) & $58977314$\\ \hline
5&(10,2,4) & $487020090$\\ \hline
5&(12,1,4) & $3986631790$\\ \hline
5&(14,0,4) & $32664263244$\\ \hline
5&(1,8,3) & $168160$\\ \hline
5&(3,7,3) & $1426788$\\ \hline
5&(5,6,3) & $12341640$\\ \hline
5&(7,5,3) & $106742892$\\ \hline
5&(9,4,3) & $917273760$\\ \hline
5&(11,3,3) & $7835510640$\\ \hline
5&(13,2,3) & $66838183448$\\ \hline
5&(15,1,3) & $572970400800$\\ \hline
5&(17,0,3) & $4963870717184$\\ \hline
\end{tabular}
\end{table}
\begin{table}[H]
\centering
\caption{$g=0,M_6^1,d=5$ part 2} \label{g=0,N=6,k=1,d=5,2}
\begin{tabular}{|c|l|l|l|l|l|}
\hline
d&(a,b,c)& $N^0_{d,a,b,c}$  \\
\hline
5&(0,10,2) & $275340$\\ \hline
5&(2,9,2) & $2427884$\\ \hline
5&(4,8,2) & $21445040$\\ \hline
5&(6,7,2) & $190810312$\\ \hline
5&(8,6,2) & $1697371800$\\ \hline
5&(10,5,2) & $15059634800$\\ \hline
5&(12,4,2) & $133462672144$\\ \hline
5&(14,3,2) & $1185922233290$\\ \hline
5&(16,2,2) & $10614063989964$\\ \hline
5&(18,1,2) & $96073499325220$\\ \hline
5&(20,0,2) & $882272821107200$\\ \hline
5&(1,11,1) & $3941780$\\ \hline
5&(3,10,1) & $36486980$\\ \hline
5&(5,9,1) & $335462284$\\ \hline
5&(7,8,1) & $3085793380$\\ \hline
5&(9,7,1) & $28364597480$\\ \hline
5&(11,6,1) & $260604570680$\\ \hline
5&(13,5,1) & $2398306990560$\\ \hline
5&(15,4,1) & $22179879220568$\\ \hline
5&(17,3,1) & $206859628175260$\\ \hline
5&(19,2,1) & $1951736958419580$\\ \hline
5&(21,1,1) & $18676573063528460$\\ \hline
5&(23,0,1) & $181610832693333060$\\ \hline
5&(0,13,0) & $5953000$\\ \hline
5&(2,12,0) & $59294040$\\ \hline
5&(4,11,0) & $573878820$\\ \hline
5&(6,10,0) & $5489009900$\\ \hline
5&(8,9,0) & $52316386080$\\ \hline
5&(10,8,0) & $497999093480$\\ \hline
5&(12,7,0) & $4743580571280$\\ \hline
5&(14,6,0) & $45328575942720$\\ \hline
5&(16,5,0) & $435813009759000$\\ \hline
5&(18,4,0) & $4228201646521080$\\ \hline
5&(20,3,0) & $41500751630424420$\\ \hline
5&(22,2,0) & $412968880889100580$\\ \hline
5&(24,1,0) & $4173087641902059600$\\ \hline
5&(26,0,0) & $42876778851631702000$\\ \hline
\end{tabular}
\end{table}
\begin{table}[H]
\centering
\caption{$g=1,M_6^1,d=1,2$} \label{g=1,N=6,k=1}
\begin{tabular}{|c|l|l|l|l|l|}
\hline
d&(a,b,c)& $N^1_{d,a,b,c}$& $E_{d,a,b,c}$\\
\hline
1&(0,1,1) & $-\frac{5 }{24}$&$0$\\ \hline
1&(2,0,1) & $-\frac{5 }{24}$&$0$\\ \hline
1&(1,2,0) & $-\frac{5}{12}$&$0$\\ \hline
1&(3,1,0) & $-\frac{5 }{8}$&$0$\\ \hline
1&(5,0,0) & $-\frac{25}{24}$&$0$\\ \hline
2&(1,0,3) & $-\frac{3 }{8}$&$0$\\ \hline
2&(0,2,2) & $-\frac{3}{8}  $&$0$\\ \hline
2&(2,1,2) & $-\frac{17}{12} $&$0$\\ \hline
2&(4,0,2) & $-\frac{91}{24} $&$0$\\ \hline
2&(1,3,1) & $-\frac{25}{12} $&$0$\\ \hline
2&(3,2,1) & $-\frac{173}{24} $&$0$\\ \hline
2&(5,1,1) & $-\frac{545}{24} $&$0$\\ \hline
2&(7,0,1) & $-\frac{1765 }{24}$&$0$\\ \hline
2&(0,5,0) & $-\frac{10}{3} $&$0$\\ \hline
2&(2,4,0) & $-\frac{73}{6} $&$0$\\ \hline
2&(4,3,0) & $-\frac{89}{2}$&$0$\\ \hline
2&(6,2,0) & $-\frac{3805}{24} $&$0$\\ \hline
2&(8,1,0) & $-\frac{1735}{3}$&$0$\\ \hline
2&(10,0,0) & $-2200$&0\\ \hline
\end{tabular}
\end{table}

\begin{table}[H]
\centering
\caption{$g=1,M_6^1,d=3$} \label{g=1,N=6,k=1,d=3}
\begin{tabular}{|c|l|l|l|l|l|}
\hline
d&(a,b,c)& $N^1_{d,a,b,c}$ &$E_{d,a,b,c}$ \\
\hline
3&(0,0,5) & $0 $&$0$\\ \hline
3&(1,1,4) & $-\frac{55}{24} $&$0$\\ \hline
3&(3,0,4) & $-\frac{55}{4} $&$0$\\ \hline
3&(0,3,3) & $-\frac{95}{24} $&$0$\\ \hline
3&(2,2,3) & $-\frac{475}{24}  $&$0$\\ \hline
3&(4,1,3) & $-\frac{785}{8}  $&$0$\\ \hline
3&(6,0,3) & $-\frac{10337}{24}$&$1$\\ \hline
3&(1,4,2) & $-\frac{65}{2} $&$0$\\ \hline
3&(3,3,2) & $-\frac{3935}{24}  $&$0$\\ \hline
3&(5,2,2) & $-\frac{6519}{8} $&$1$\\ \hline
3&(7,1,2) & $-\frac{15603}{4}  $&$14$\\ \hline
3&(9,0,2) & $-\frac{75273}{4}  $&$114$\\ \hline
3&(0,6,1) & $-\frac{160}{3}  $&$0$\\ \hline
3&(2,5,1) & $-\frac{830}{3}  $&$0$\\ \hline
3&(4,4,1) & $-\frac{4373}{3}  $&$2$\\ \hline
3&(6,3,1) & $-\frac{182575}{24}  $&$25$\\ \hline
3&(8,2,1) & $-\frac{236081}{6}   $&$222$\\ \hline
3&(10,1,1) & $-\frac{1236305}{6}$&$1650$\\ \hline
3&(12,0,1) & $-\frac{2212985 }{2} $&$11325$\\ \hline
3&(1,7,0) & $-\frac{2705 }{6}  $&$0$\\ \hline
3&(3,6,0) & $-\frac{7600}{3} $&$5$\\ \hline
3&(5,5,0) & $-\frac{84815}{6}   $&$55$\\ \hline
3&(7,4,0) & $-\frac{235232}{3} $&$468$\\ \hline
3&(9,3,0) & $-\frac{2613469}{6}  $&$3558$\\ \hline
3&(11,2,0) & $-\frac{29454205}{12}  $&$25275$\\ \hline
3&(13,1,0) & $-\frac{56619575 }{4}$&$173490$\\ \hline
3&(15,0,0) & $-\frac{503847475 }{6} $&$1175300$\\ \hline
\end{tabular}
\end{table}
\end{multicols}
\begin{table}[H]
\centering
\caption{$g=1,M_6^1,d=4$} \label{g=1,N=6,k=1,d=4}
\begin{tabular}{|c|l|l|l|l|l|}
\hline
d&(a,b,c)& $N^1_{d,a,b,c}$ &$E_{d,a,b,c}$ \\
\hline
4&(0,1,6) & $-\frac{35}{8}$&$0$\\ \hline
4&(2,0,6) & $-\frac{725}{24}$&$0$\\ \hline
4&(1,2,5) & $-\frac{1255}{24}$&$0$\\ \hline
4&(3,1,5) & $-\frac{2165}{6}$&$0$\\ \hline
4&(5,0,5) & $-\frac{28625}{12}$&$0$\\ \hline
4&(0,4,4) & $-\frac{341}{4}$&$1$\\ \hline
4&(2,3,4) & $-\frac{6991}{12}$&$4$\\ \hline
4&(4,2,4) & $-\frac{31035}{8}$&$32$\\ \hline
4&(6,1,4) & $-\frac{594955}{24}$&$310$\\ \hline
4&(8,0,4) & $-\frac{449893}{3}$&$3220$\\ \hline
4&(1,5,3) & $-\frac{2881}{3}$&$14$\\ \hline
4&(3,4,3) & $-\frac{25957}{4}$&$96$\\ \hline
4&(5,3,3) & $-\frac{173685}{4}$&$785$\\ \hline
4&(7,2,3) & $-\frac{1135543}{4}$&$6755$\\ \hline
4&(9,1,3) & $-1819093$&$57960$\\ \hline
4&(11,0,3) & $-\frac{70041271}{6}$&$ 473586$\\ \hline
4&(0,7,2) & $-\frac{9647}{6} $&$29$\\ \hline
4&(2,6,2) & $-10889 $&$228$\\ \hline
4&(4,5,2) & $-75377$&$1799$\\ \hline
4&(6,4,2) & $-\frac{6227887}{12}$&$14745$\\ \hline
4&(8,3,2) & $-\frac{10604846}{3}$&$122155$\\ \hline
4&(10,2,2) & $-\frac{287430247}{12}$&$1004916$\\ \hline
4&(12,1,2) & $-\frac{1961004899}{12}$&$8127687$\\ \hline
4&(14,0,2) & $-\frac{13594466647}{12}$&$65017656$\\ \hline
4&(1,8,1) & $-\frac{107875}{6}$&$460$\\ \hline
4&(3,7,1) & $-\frac{386698}{3}$&$3837$\\ \hline
4&(5,6,1) & $-\frac{2782210}{3}$&$31470$\\ \hline
4&(7,5,1) & $-\frac{19937800}{3} $&$259325$\\ \hline
4&(9,4,1) & $-47440464 $&$2139212$\\ \hline
4&(11,3,1) & $-\frac{1017676760}{3}$&$ 17589245$\\ \hline
4&(13,2,1) & $-\frac{9795944647}{4} $&$144007483$\\ \hline
4&(15,1,1) & $-\frac{53891135695}{3} $&$ 1178586170$\\ \hline
4&(17,0,1) & $-\frac{806789380135}{6} $&$9691893740$\\ \hline
4&(0,10,0) & $-\frac{164635}{6}$&$870$\\ \hline
4&(2,9,0) & $-\frac{639692}{3}$&$7708$\\ \hline
4&(4,8,0) & $-\frac{9683117}{6}$&$65124$\\ \hline
4&(6,7,0) & $-\frac{72844595}{6}$&$ 543660$\\ \hline
4&(8,6,0) & $-\frac{272974555}{3}$&$4535180$\\ \hline
4&(10,5,0) & $-\frac{2045100320}{3}$&$37854490$\\ \hline
4&(12,4,0) & $-\frac{10265421909}{2}$&$316190712$\\ \hline
4&(14,3,0) & $-\frac{78131149107}{2}$&$2646562486$\\ \hline
4&(16,2,0) & $-\frac{604014726315}{2}$&$22264309750$\\ \hline
4&(18,1,0) & $-2379763852305 $&$188886527100$\\ \hline
4&(20,0,0) & $-\frac{114958360544525}{6}$&$1620988570200$\\ \hline
\end{tabular}
\end{table}
\begin{table}[H]
\centering
\caption{$g=1,M_6^1,d=5$, part 1} \label{g=1,N=6,k=1,d=5,part1}
\begin{tabular}{|c|l|l|l|l|l|}
\hline
d&(a,b,c)& $N^1_{d,a,b,c}$ &$E_{d,a,b,c}$ \\
\hline
5&(1,0,8) & $-\frac{1945}{24}$&$3$\\ \hline
5&(0,2,7) & $-\frac{1435}{12}$&$8$\\ \hline
5&(2,1,7) & $-\frac{2405}{2}$&$45$\\ \hline
5&(4,0,7) & $-\frac{126757}{12}$&$324$\\ \hline
5&(1,3,6) & $-\frac{15315}{8}$&$105$\\ \hline
5&(3,2,6) & $-\frac{50707}{3}$&$771$\\ \hline
5&(5,1,6) & $-\frac{3386155}{24}$&$6330$\\ \hline
5&(7,0,6) & $-\frac{2231555}{2}$&$55935$\\ \hline
5&(0,5,5) & $-3205$&$ 220$\\ \hline
5&(2,4,5) & $-\frac{54877}{2}$&$1694$\\ \hline
5&(4,3,5) & $-\frac{5573975}{24}$&$14090$\\ \hline
5&(6,2,5) & $-\frac{15246635}{8}$&$123855$\\ \hline
5&(8,1,5) & $-\frac{181156655}{12}$&$ 1127945$\\ \hline
5&(10,0,5) & $-\frac{1390208755}{12}$&$10460910$\\ \hline
5&(1,6,4) & $-\frac{274369}{6}$&$3401$\\ \hline
5&(3,5,4) & $-\frac{1160363}{3}$&$29084$\\ \hline
5&(5,4,4) & $-\frac{6522287}{2}$&$ 257413$\\ \hline
5&(7,3,4) & $-\frac{324143825}{12}$&$2335035$\\ \hline
5&(9,2,4) & $-\frac{2632411415}{12}$&$21477895$\\ \hline
5&(11,1,4) & $-\frac{7029536823}{4}$&$ 198302628$\\ \hline
5&(13,0,4) & $-\frac{84374784907}{6}$&$1821344223$\\ \hline
5&(0,8,3) & $-77296$&$6224$\\ \hline
5&(2,7,3) & $-\frac{1297759}{2}$&$56623$\\ \hline
5&(4,6,3) & $-5569635 $&$513600$\\ \hline
5&(6,5,3) & $-\frac{143175565}{3}$&$4719945$\\ \hline
5&(8,4,3) & $-\frac{1214792630}{3}$&$43829370$\\ \hline
5&(10,3,3) & $-\frac{6804949031}{2}$&$409413777$\\ \hline
5&(12,2,3) & $-\frac{85400501522}{3}$&$ 3832247481$\\ \hline
5&(14,1,3) & $-\frac{1434207179585}{6}$&$ 35880218130$\\ \hline
5&(16,0,3) & $-2028186128080$&$336649015020$\\ \hline
5&(1,9,2) & $-\frac{2158487}{2}$&$104214$\\ \hline
5&(3,8,2) & $-\frac{28335680}{3}$&$986530$\\ \hline
5&(5,7,2) & $-83306285$&$9296155$\\ \hline
5&(7,6,2) & $-\frac{2200182005}{3}$&$87964690$\\ \hline
5&(9,5,2) & $-6423114302$&$836194988$\\ \hline
5&(11,4,2) & $-56047183707$&$7977161338$\\ \hline
5&(13,3,2) & $-\frac{5873452441955}{12}$&$76321377460$\\ \hline
5&(15,2,2) & $-4301607795510$&$732727776780$\\ \hline
5&(17,1,2) & $-\frac{229390883916695}{6}$&$7072862332690$\\ \hline
5&(19,0,2) & $-\frac{1034728934410532}{3}$&$68819539583286$\\ \hline
\end{tabular}
\end{table}
\begin{table}[H]
\centering
\caption{$g=1,M_6^1,d=5$, part 2} \label{g=1,N=6,k=1,d=5,part2}
\begin{tabular}{|c|l|l|l|l|l|}
\hline
d&(a,b,c)& $N^1_{d,a,b,c}$ &$E_{d,a,b,c}$ \\
\hline
5&(0,11,1) & $-\frac{3402365}{2}$&$184480$\\ \hline
5&(2,10,1) & $-\frac{94016375}{6}$&$1825740$\\ \hline
5&(4,9,1) & $-\frac{857749505}{6}$&$17812920$\\ \hline
5&(6,8,1) & $-\frac{7821853175}{6}$&$172988900$\\ \hline
5&(8,7,1) & $-11859338793$&$1681583302$\\ \hline
5&(10,6,1) & $-\frac{322851405685}{3}$&$16388989480$\\ \hline
5&(12,5,1) & $-\frac{2929215576475}{3}$&$160253018950$\\ \hline
5&(14,4,1) & $-8891388713520$&$1573356927950$\\ \hline
5&(16,3,1) & $-\frac{489643030178765}{6}$&$15531173703960$\\ \hline
5&(18,2,1) & $-\frac{4546452758023601}{6}$&$154438870292624$\\ \hline
5&(20,1,1) & $-\frac{14277641483246075}{2}$&$1550155293988500$\\ \hline
5&(22,0,1) & $-\frac{410343935038064375}{6}$&$15736008036054900$\\ \hline
5&(1,12,0) & $-\frac{73969625}{3}$&$3274890$\\ \hline
5&(3,11,0) & $-\frac{1428317155}{6}$&$33189090$\\ \hline
5&(5,10,0) & $-\frac{13579189585}{6}$&$332653440$\\ \hline
5&(7,9,0) & $-\frac{64202869786}{3}$&$3320681718$\\ \hline
5&(9,8,0) & $-201753545535$&$33152681950$\\ \hline
5&(11,7,0) & $-1900031299700$&$331728574700$\\ \hline
5&(13,6,0) & $-17926222775700$&$3331531685250$\\ \hline
5&(15,5,0) & $-\frac{510047851228460}{3}$&$33628906542030$\\ \hline
5&(17,4,0) & $-\frac{4879475383947517}{3}$&$341728487604966$\\ \hline
5&(19,3,0) & $-\frac{31485553789132145}{2}$&$3501911868774390$\\ \hline
5&(21,2,0) & $-\frac{927227199825166475}{6}$&$ 36252054632249100$\\ \hline
5&(23,1,0) & $-\frac{4624335742996347470}{3}$&$379702695213071370$\\ \hline
5&(25,0,0) & $-15644867027685574170 $&$4029226126511838330$\\ \hline
\end{tabular}
\end{table}
\newpage
\begin{multicols}{2}
\begin{table}[H]
\centering
\caption{$g=0,M_6^2,d=1,2,3$} \label{g=0,N=6,k=2,d=1,2,3}
\begin{tabular}{|c|l|l|l|l|l|}
\hline
d&(a,b,c)& $N^0_{d,a,b,c}$  \\
\hline
1&(0,1,1) & $4$\\ \hline
1&(2,0,1) & $4$\\ \hline
1&(1,2,0) & $8$\\ \hline
1&(3,1,0) & $12$\\ \hline
1&(5,0,0) & $20$\\ \hline
2&(0,0,3) & $8$\\ \hline
2&(1,1,2) & $16$\\ \hline
2&(3,0,2) & $32$\\ \hline
2&(0,3,1) & $16$\\ \hline
2&(2,2,1) & $64$\\ \hline
2&(4,1,1) & $184$\\ \hline
2&(6,0,1) & $576$\\ \hline
2&(1,4,0) & $96$\\ \hline
2&(3,3,0) & $368$\\ \hline
2&(5,2,0) & $1280$\\ \hline
2&(7,1,0) & $4632$\\ \hline
2&(9,0,0) & $17704$\\ \hline
3&(1,0,4) & $64$\\ \hline
3&(0,2,3) & $64$\\ \hline
3&(2,1,3) & $320$\\ \hline
3&(4,0,3) & $1152$\\ \hline
3&(1,3,2) & $448$\\ \hline
3&(3,2,2) & $2240$\\ \hline
3&(5,1,2) & $9888$\\ \hline
3&(7,0,2) & $45504$\\ \hline
3&(0,5,1) & $640$\\ \hline
3&(2,4,1) & $3712$\\ \hline
3&(4,3,1) & $19552$\\ \hline
3&(6,2,1) & $98272$\\ \hline
3&(8,1,1) & $505552$\\ \hline
3&(10,0,1) & $2702000$\\ \hline
3&(1,6,0) & $6144$\\ \hline
3&(3,5,0) & $35584$\\ \hline
3&(5,4,0) & $199424$\\ \hline
3&(7,3,0) & $1102752$\\ \hline
3&(9,2,0) & $6213728$\\ \hline
3&(11,1,0) & $36112336$\\ \hline
3&(13,0,0) & $217541136$\\ \hline
\end{tabular}
\end{table}

\begin{table}[H]
\centering
\caption{$g=0,M_6^2, d=4$} \label{g=0,N=6,k=2,d=4}
\begin{tabular}{|c|l|l|l|l|l|}
\hline
d&(a,b,c)& $N^0_{d,a,b,c}$  \\
\hline
4&(0,1,5) & $384$\\ \hline
4&(2,0,5) & $2560$\\ \hline
4&(1,2,4) & $3328$\\ \hline
4&(3,1,4) & $20096$\\ \hline
4&(5,0,4) & $103296$\\ \hline
4&(0,4,3) & $4608$\\ \hline
4&(2,3,3) & $31616$\\ \hline
4&(4,2,3) & $198912$\\ \hline
4&(6,1,3) & $1177152$\\ \hline
4&(8,0,3) & $7139008$\\ \hline
4&(1,5,2) & $48640$\\ \hline
4&(3,4,2) & $345088$\\ \hline
4&(5,3,2) & $2312960$\\ \hline
4&(7,2,2) & $15125760$\\ \hline
4&(9,1,2) & $100512384$\\ \hline
4&(11,0,2) & $687547904$\\ \hline
4&(0,7,1) & $74240$\\ \hline
4&(2,6,1) & $579584$\\ \hline
4&(4,5,1) & $4263424$\\ \hline
4&(6,4,1) & $30403072$\\ \hline
4&(8,3,1) & $215004288$\\ \hline
4&(10,2,1) & $1541202944$\\ \hline
4&(12,1,1) & $11310265600$\\ \hline
4&(14,0,1) & $85345355904$\\ \hline
4&(1,8,0) & $980992$\\ \hline
4&(3,7,0) & $7649280$\\ \hline
4&(5,6,0) & $58502144$\\ \hline
4&(7,5,0) & $441611520$\\ \hline
4&(9,4,0) & $3335082240$\\ \hline
4&(11,3,0) & $25524363520$\\ \hline
4&(13,2,0) & $199433493504$\\ \hline
4&(15,1,0) & $1597307381568$\\ \hline
4&(17,0,0) & $13143513966080$\\ \hline
\end{tabular}
\end{table}

\begin{table}[H]
\centering
\caption{$g=0,M_6^2, d=5$, part 1} \label{g=0,N=6,k=2,d=5, part1}
\begin{tabular}{|c|l|l|l|l|l|}
\hline
d&(a,b,c)& $N^0_{d,a,b,c}$  \\
\hline
5&(0,0,7) & $3072$\\ \hline
5&(1,1,6) & $33280$\\ \hline
5&(3,0,6) & $238592$\\ \hline
5&(0,3,5) & $45056$\\ \hline
5&(2,2,5) & $357376$\\ \hline
5&(4,1,5) & $2592000$\\ \hline
5&(6,0,5) & $17329920$\\ \hline
5&(1,4,4) & $530432$\\ \hline
5&(3,3,4) & $4331520$\\ \hline
5&(5,2,4) & $33241088$\\ \hline
5&(7,1,4) & $246732544$\\ \hline
5&(9,0,4) & $1859635968$\\ \hline
5&(0,6,3) & $774144$\\ \hline
5&(2,5,3) & $6976512$\\ \hline
5&(4,4,3) & $59192320$\\ \hline
5&(6,3,3) & $482989056$\\ \hline
5&(8,2,3) & $3885555456$\\ \hline
5&(10,1,3) & $31625447936$\\ \hline
5&(12,0,3) & $263294153728$\\ \hline
5&(1,7,2) & $11044864$\\ \hline
5&(3,6,2) & $101627904$\\ \hline
5&(5,5,2) & $898333696$\\ \hline
5&(7,4,2) & $7760710656$\\ \hline
5&(9,3,2) & $66707246848$\\ \hline
5&(11,2,2) & $579410893568$\\ \hline
5&(13,1,2) & $5127561935232$\\ \hline
5&(15,0,2) & $46417851808512$\\ \hline
5&(0,9,1) & $17539072$\\ \hline
5&(2,8,1) & $171974656$\\ \hline
5&(4,7,1) & $1620875264$\\ \hline
5&(6,6,1) & $14922719232$\\ \hline
5&(8,5,1) & $135757207552$\\ \hline
5&(10,4,1) & $1235286731776$\\ \hline
5&(12,3,1) & $11356463308416$\\ \hline
5&(14,2,1) & $106142531433856$\\ \hline
5&(16,1,1) & $1012183022176064$\\ \hline
5&(18,0,1) & $9868286056471104$\\ \hline
\end{tabular}
\end{table}
\begin{table}[H]
\centering
\caption{$g=0,M_6^2, d=5$, part 2} \label{g=0,N=6,k=2,d=5, part2}
\begin{tabular}{|c|l|l|l|l|l|}
\hline
d&(a,b,c)& $N^0_{d,a,b,c}$  \\
\hline
5&(1,10,0) & $295157760$\\ \hline
5&(3,9,0) & $2892357632$\\ \hline
5&(5,8,0) & $28004556800$\\ \hline
5&(7,7,0) & $268191540224$\\ \hline
5&(9,6,0) & $2557415162880$\\ \hline
5&(11,5,0) & $24477957189120$\\ \hline
5&(13,4,0) & $236783110502912$\\ \hline
5&(15,3,0) & $2325894666007680$\\ \hline
5&(17,2,0) & $23271693753656448$\\ \hline
5&(19,1,0) & $237640853649062080$\\ \hline
5&(21,0,0) & $2479850169230342720$\\ \hline
\end{tabular}
\end{table}
\begin{table}[H]
\centering
\caption{$g=1,M_6^2, d=1,2,3$} \label{g=1,N=6,k=2,d=1,2,3}
\begin{tabular}{|c|l|l|l|l|l|}
\hline
d&(a,b,c)& $N^1_{d,a,b,c}$&$E_{d,a,b,c}$  \\
\hline
1&(1,0,1) & $-\frac{1}{2}$&$0$\\ \hline
1&(0,2,0) & $-1$&$0$\\ \hline
1&(2,1,0) & $-\frac{3 }{2}$&$0$\\ \hline
1&(4,0,0) & $-\frac{5 }{2}$&$0$\\ \hline
2&(0,1,2) & $-\frac{10}{3} $&$0$\\ \hline
2&(2,0,2) & $-\frac{20}{3}  $&$0$\\ \hline
2&(1,2,1) & $-\frac{40}{3} $&$0$\\ \hline
2&(3,1,1) & $-\frac{115}{3} $&$0$\\ \hline
2&(5,0,1) & $-120  $&$0$\\ \hline
2&(0,4,0) & $-20 $&$0$\\ \hline
2&(2,3,0) & $-\frac{230}{3}  $&$0$\\ \hline
2&(4,2,0) & $-\frac{800}{3}  $&$0$\\ \hline
2&(6,1,0) & $-965 $&$0$\\ \hline
2&(8,0,0) & $-\frac{11065 }{3}$&$0$\\ \hline
3&(0,0,4) & $-\frac{56 }{3}$&$0$\\ \hline
3&(1,1,3) & $-88 $&$0$\\ \hline
3&(3,0,3) & $-\frac{928}{3}$&$0$\\ \hline
3&(0,3,2) & $-120$&$0$\\ \hline
3&(2,2,2) & $-600$&$0$\\ \hline
3&(4,1,2) & $-\frac{7852}{3}$&$0$\\ \hline
3&(6,0,2) & $-11936$&$0$\\ \hline
3&(1,4,1) & $-\frac{2944}{3} $&$0$\\ \hline
3&(3,3,1) & $-\frac{15476}{3}$&$0$\\ \hline
3&(5,2,1) & $-25756$&$0$\\ \hline
3&(7,1,1) & $-\frac{395158}{3}$&$8$\\ \hline
3&(9,0,1) & $-\frac{2101942 }{3}$&$136$\\ \hline
3&(0,6,0) & $-\frac{4832}{3}$&$0$\\ \hline
3&(2,5,0) & $-\frac{27952}{3}$&$0$\\ \hline
3&(4,4,0) & $-52112$&$0$\\ \hline
3&(6,3,0) & $-\frac{861020}{3}$&$16$\\ \hline
3&(8,2,0) & $-1611044$&$272$\\ \hline
3&(10,1,0) & $-9333366$&$3000$\\ \hline
3&(12,0,0) & $-\frac{168264818}{3}$&$27240$\\ \hline
\end{tabular}
\end{table}

\begin{table}[H]
\centering
\caption{$g=1,M_6^2, d=4$} \label{g=1,N=6,k=2,d=4}
\begin{tabular}{|c|l|l|l|l|l|}
\hline
d&(a,b,c)& $N^1_{d,a,b,c}$&$E_{d,a,b,c}$  \\
\hline
4&(1,0,5) & $-\frac{2560}{3}$&$0$\\ \hline
4&(0,2,4) & $-1056$&$0$\\ \hline
4&(2,1,4) & $-6320$&$0$\\ \hline
4&(4,0,4) & $-\frac{94192}{3}$&$256$\\ \hline
4&(1,3,3) & $-9680 $&$0$\\ \hline
4&(3,2,3) & $-\frac{181024}{3}$&$256$\\ \hline
4&(5,1,3) & $-\frac{1050152}{3}$&$3584$\\ \hline
4&(7,0,3) & $-2088664$&$30720$\\ \hline
4&(0,5,2) & $-\frac{43520}{3}$&$0$\\ \hline
4&(2,4,2) & $-\frac{308864}{3}$&$512$\\ \hline
4&(4,3,2) & $-\frac{2053856}{3}$&$6144$\\ \hline
4&(6,2,2) & $-\frac{13260032}{3}$&$59392$\\ \hline
4&(8,1,2) & $-\frac{87089248}{3} $&$499296$\\ \hline
4&(10,0,2) & $-\frac{589841888}{3} $&$4007808$\\ \hline
4&(1,6,1) & $-169856$&$1280$\\ \hline
4&(3,5,1) & $-\frac{3742112}{3}$&$12800$\\ \hline
4&(5,4,1) & $-8838784$&$ 118784$\\ \hline
4&(7,3,1) & $-\frac{185860784}{3}$&$ 1043936$\\ \hline
4&(9,2,1) & $-440296032 $&$8781120$\\ \hline
4&(11,1,1) & $-3206552104$&$72704112$\\ \hline
4&(13,0,1) & $-\frac{72109405664}{3}$&$602772128$\\ \hline
4&(0,8,0) & $-283776$&$3584$\\ \hline
4&(2,7,0) & $-2210240$&$28672$\\ \hline
4&(4,6,0) & $-16860672$&$251904$\\ \hline
4&(6,5,0) & $-\frac{379892480}{3}$&$2211200$\\ \hline
4&(8,4,0) & $-\frac{2850607936}{3}$&$19140736$\\ \hline
4&(10,3,0) & $-7225953840$&$163776288$\\ \hline
4&(12,2,0) & $-\frac{168399937952}{3}$&$1399658624$\\ \hline
4&(14,1,0) & $-\frac{1341957809680}{3}$&$12047200912$\\ \hline
4&(16,0,0) & $-3665151541000$&$105004257360$\\ \hline
\end{tabular}
\end{table}
\end{multicols}
\begin{table}[H]
\centering
\caption{$g=1,M_6^2, d=5$} \label{g=1,N=6,k=2,d=5}
\begin{tabular}{|c|l|l|l|l|l|}
\hline
d&(a,b,c)& $N^1_{d,a,b,c}$&$E_{d,a,b,c}$  \\
\hline
5&(0,1,6) & $-11840$&$0$\\ \hline
5&(2,0,6) & $-84224$&$0$\\ \hline
5&(1,2,5) & $-120832$&$768$\\ \hline
5&(3,1,5) & $-\frac{2583392}{3}$&$13824$\\ \hline
5&(5,0,5) & $-\frac{16612064}{3}$&$192512$\\ \hline
5&(0,4,4) & $-\frac{518144}{3}$&$2048$\\ \hline
5&(2,3,4) & $-\frac{4207744}{3}$&$28416$\\ \hline
5&(4,2,4) & $-\frac{31740928}{3}$&$328960$\\ \hline
5&(6,1,4) & $-\frac{229497440}{3}$&$3457152$\\ \hline
5&(8,0,4) & $-\frac{1689720416}{3}$&$32359808$\\ \hline
5&(1,5,3) & $-2195328$&$ 60928$\\ \hline
5&(3,4,3) & $-18503296$&$649728$\\ \hline
5&(5,3,3) & $-\frac{446421376}{3}$&$6655104$\\ \hline
5&(7,2,3) & $-1174715232$&$65456768$\\ \hline
5&(9,1,3) & $-9391653184$&$618276096$\\ \hline
5&(11,0,3) & $-76961465920$&$5753289984$\\ \hline
5&(0,7,2) & $-3375616 $&$132096$\\ \hline
5&(2,6,2) & $-31075840$&$1330176$\\ \hline
5&(4,5,2) & $-272878080$&$13367296$\\ \hline
5&(6,4,2) & $-\frac{6989383808}{3}$&$132810240$\\ \hline
5&(8,3,2) & $-\frac{59226474208}{3}$&$1299074944$\\ \hline
5&(10,2,2) & $-169095332896$&$12544392576$\\ \hline
5&(12,1,2) & $-\frac{4432072954096}{3}$&$120900975936$\\ \hline
5&(14,0,2) & $-\frac{39663151820800}{3}$&$1173261629056$\\ \hline
5&(1,8,1) & $-\frac{154276864}{3}$&$2738176$\\ \hline
5&(3,7,1) & $-\frac{1451733632}{3}$&$27183616$\\ \hline
5&(5,6,1) & $-\frac{13282234496}{3}$&$ 271901184$\\ \hline
5&(7,5,1) & $-\frac{119682815168}{3} $&$2715012608$\\ \hline
5&(9,4,1) & $-\frac{1076978890048}{3} $&$27008532736$\\ \hline
5&(11,3,1) & $-\frac{9791608439120}{3} $&$ 268391724096$\\ \hline
5&(13,2,1) & $-\frac{90566832590384}{3} $&$2679621302336$\\ \hline
5&(15,1,1) & $-\frac{855469605426472}{3} $&$27017537730048$\\ \hline
5&(17,0,1) & $-\frac{8269323668574136}{3} $&$276138960259456$\\ \hline
5&(0,10,0) & $-86736896$&$5566464$\\ \hline
5&(2,9,0) & $-\frac{2547003392}{3}$&$55035904$\\ \hline
5&(4,8,0) & $-\frac{24584536064}{3}$&$555792384$\\ \hline
5&(6,7,0) & $-78025121792$&$5650489344$\\ \hline
5&(8,6,0) & $-\frac{2214589877504}{3}$&$57577844736$\\ \hline
5&(10,5,0) & $-7003208209024$&$587949802496$\\ \hline
5&(12,4,0) & $-67142653348736$&$6030797709312$\\ \hline
5&(14,3,0) & $-\frac{1961908485684176}{3}$&$62362747054848$\\ \hline
5&(16,2,0) & $-6492457837633744$&$652295164595968$\\ \hline
5&(18,1,0) & $-65834214232016696$&$6919499749525632$\\ \hline
5&(20,0,0) & $-682743285735788648 $&$74585887608233472$\\ \hline
\end{tabular}
\end{table}
\begin{multicols}{2}
\begin{table}[H]
\centering
\caption{$g=0,M_6^3,d=1,2,3$} \label{g=0, N=6,k=3,d=1,2,3}
\begin{tabular}{|c|l|l|l|l|l|}
\hline
d&(a,b,c)& $N^0_{d,a,b,c}$  \\
\hline
1&(1,0,1) & $18$\\ \hline
1&(0,2,0) & $45$\\ \hline
1&(2,1,0) & $63$\\ \hline
1&(4,0,0) & $108$\\ \hline
2&(1,0,2) & $108$\\ \hline
2&(0,2,1)&$378$ \\ \hline
2&(2,1,1) & $864$\\ \hline
2&(4,0,1) & $2754$\\ \hline
2&(1,3,0) & $2187$\\ \hline
2&(3,2,0) & $7047$\\ \hline
2&(5,1,0) & $25758$\\ \hline
2&(7,0,0) & $102060$\\ \hline
3&(1,0,3) & $1944$\\ \hline
3&(0,2,2) & $7452$\\ \hline
3&(2,1,2) & $24300$\\ \hline
3&(4,0,2) & $110808$\\ \hline
3&(1,3,1) & $65610$\\ \hline
3&(3,2,1) & $290142$\\ \hline
3&(5,1,1) & $1469664$\\ \hline
3&(7,0,1) & $8019000$\\ \hline
3&(0,5,0) & $131220$\\ \hline
3&(2,4,0) & $715149$\\ \hline
3&(4,3,0) & $3766014$\\ \hline
3&(6,2,0) & $21210984$\\ \hline
3&(8,1,0) & $126574812$\\ \hline
3&(10,0,0) & $796767840$\\ \hline
\end{tabular}
\end{table}

\begin{table}[H]
\centering
\caption{$g=0,M_6^3,d=4$} \label{g=0, N=6,k=3,d=4}
\begin{tabular}{|c|l|l|l|l|l|}
\hline
d&(a,b,c)& $N^0_{d,a,b,c}$  \\
\hline
4&(1,0,4) & $62208$\\ \hline
4&(0,2,3) & $248832$\\ \hline
4&(2,1,3) & $1057536$\\ \hline
4&(4,0,3) & $6216912$\\ \hline
4&(1,3,2) & $2991816$\\ \hline
4&(3,2,2) & $16726176$\\ \hline
4&(5,1,2) & $107372952$\\ \hline
4&(7,0,2) & $741917880$\\ \hline
4&(0,5,1) & $6495390$\\ \hline
4&(2,4,1) & $42554646$\\ \hline
4&(4,3,1) & $278396352$\\ \hline
4&(6,2,1) & $1955545416$\\ \hline
4&(8,1,1) & $14516081280$\\ \hline
4&(10,0,1) & $112889878200$\\ \hline
4&(1,6,0) & $94242204$\\ \hline
4&(3,5,0) & $683570907$\\ \hline
4&(5,4,0) & $4991884362$\\ \hline
4&(7,3,0) & $38081408688$\\ \hline
4&(9,2,0) & $303386439924$\\ \hline
4&(11,1,0) & $2518696831248$\\ \hline
4&(13,0,0) & $21763790077104$\\ \hline
\end{tabular}
\end{table}

\begin{table}[H]
\centering
\caption{$g=0,M_6^3,d=5$} \label{g=0, N=6,k=3,d=5}
\begin{tabular}{|c|l|l|l|l|l|}
\hline
d&(a,b,c)& $N^0_{d,a,b,c}$  \\
\hline
5&(1,0,5) & $2916000$\\ \hline
5&(0,2,4) & $11955600$\\ \hline
5&(2,1,4) & $62694000$\\ \hline
5&(4,0,4) & $449017344$\\ \hline
5&(1,3,3) & $183760488$\\ \hline
5&(3,2,3) & $1238349384$\\ \hline
5&(5,1,3) & $9579584880$\\ \hline
5&(7,0,3) & $79760239920$\\ \hline
5&(0,5,2) & $423735624$\\ \hline
5&(2,4,2) & $3238483356$\\ \hline
5&(4,3,2) & $25212033432$\\ \hline
5&(6,2,2) & $211216698360$\\ \hline
5&(8,1,2) & $1868651427024$\\ \hline
5&(10,0,2) & $17269190883936$\\ \hline
5&(1,6,1) & $7459305876$\\ \hline
5&(3,5,1) & $63063741510$\\ \hline
5&(5,4,1) & $543094438248$\\ \hline
5&(7,3,1) & $4895543320344$\\ \hline
5&(9,2,1) & $46024322567256$\\ \hline
5&(11,1,1) & $449367253012512$\\ \hline
5&(13,0,1) & $4544618405916240$\\ \hline
5&(0,8,0) & $15633399897$\\ \hline
5&(2,7,0) & $145384720047$\\ \hline
5&(4,6,0) & $1335091394610$\\ \hline
5&(6,5,0) & $12453847443000$\\ \hline
5&(8,4,0) & $119983325683644$\\ \hline
5&(10,3,0) & $1195990165100736$\\ \hline
5&(12,2,0) & $12325456056749400$\\ \hline
5&(14,1,0) & $131218747829213400$\\ \hline
5&(16,0,0) & $1442401273691663040$\\ \hline
\end{tabular}
\end{table}

\begin{table}[H]
\centering
\caption{$g=1,M_6^3, d=1,2,3$} \label{g=1,N=6,k=3, d=1,2,3}
\begin{tabular}{|c|l|l|l|l|l|}
\hline
d&(a,b,c)& $N^1_{d,a,b,c}$&$E_{d,a,b,c}$   \\
\hline
1&(0,0,1) & $-\frac{9 }{4}$&$0$\\ \hline
1&(1,1,0) & $-\frac{63 }{8}$&$0$\\ \hline
1&(3,0,0) & $-\frac{27 }{2}$&$0$\\ \hline
2&(0,0,2) & $ -18$&$0$\\ \hline
2&(1,1,1) & $-\frac{297 }{2} $&$0$\\ \hline
2&(3,0,1) & $-\frac{1917 }{4} $&$0$\\ \hline
2&(0,3,0)&$-\frac{2943}{8}$&$0$ \\ \hline
2&(2,2,0) & $-\frac{9657}{8} $&$0$\\ \hline
2&(4,1,0) & $-\frac{8883 }{2} $&$0$\\ \hline
2&(6,0,0) & $-\frac{35235 }{2}$&$0$\\ \hline
3&(0,0,3) & $-378$&$27$\\ \hline
3&(1,1,2) & $-\frac{9963}{2}$&$81$\\ \hline
3&(3,0,2) & $-22842$&$243$\\ \hline
3&(0,3,1) & $-\frac{54351}{4}$&$81$\\ \hline
3&(2,2,1) & $-\frac{239841}{4}$&$486$\\ \hline
3&(4,1,1) & $-303993$&$2187$\\ \hline
3&(6,0,1) & $-1659690$&$10935$\\ \hline
3&(1,4,0) & $-\frac{1186083}{8} $&$729$\\ \hline
3&(3,3,0) & $-\frac{3120849}{4}$&$4374$\\ \hline
3&(5,2,0) & $-4394898$&$24057$\\ \hline
3&(7,1,0) & $-\frac{52463943 }{2}$&$137781$\\ \hline
3&(9,0,0) & $-165166614$&$826686$\\ \hline
\end{tabular}
\end{table}
\end{multicols}
\begin{table}[H]
\centering
\caption{$g=1,M_6^3, d=4$} \label{g=1,N=6,k=3, d=4}
\begin{tabular}{|c|l|l|l|l|l|}
\hline
d&(a,b,c)& $N^1_{d,a,b,c}$&$E_{d,a,b,c}$   \\
\hline
4&(0,0,4) & $-12393$&$2187$\\ \hline
4&(1,1,3) & $-232065$&$13851$\\ \hline
4&(3,0,3) & $-1366875$&$67797$\\ \hline
4&(0,3,2) & $-684774$&$18954$\\ \hline
4&(2,2,2) & $-3743496$&$143613$\\ \hline
4&(4,1,2) & $-23908041$&$918540$\\ \hline
4&(6,0,2) & $-164904903$&$6259194$\\ \hline
4&(1,4,1) & $-\frac{38738331}{4}$&$253692$\\ \hline
4&(3,3,1) & $-\frac{125225433}{2}$&$1972674$\\ \hline
4&(5,2,1) & $-437451759$&$14541363$\\ \hline
4&(7,1,1) & $-3237427035$&$110257605$\\ \hline
4&(9,0,1) & $-25123505859$&$868020300$\\ \hline
4&(0,6,0) & $-\frac{42972363}{2}$&$470205$\\ \hline
4&(2,5,0) & $-\frac{1239605451}{8}$&$3965760$\\ \hline
4&(4,4,0) & $-\frac{4495376313}{4}$&$32264811$\\ \hline
4&(6,3,0) & $-\frac{17075109663}{2}$&$259323525$\\ \hline
4&(8,2,0) & $-\frac{135676682451}{2}$&$2124189360$\\ \hline
4&(10,1,0) & $-562258758393 $&$17881926768$\\ \hline
4&(12,0,0) & $-4854224373084$&$155016025290$\\ \hline
\end{tabular}
\end{table}
\begin{table}[H]
\centering
\caption{$g=1,M_6^3, d=5$} \label{g=1,N=6,k=3, d=5}
\begin{tabular}{|c|l|l|l|l|l|}
\hline
d&(a,b,c)& $N^1_{d,a,b,c}$&$E_{d,a,b,c}$   \\
\hline
5&(0,0,5) & $-538488$&$195372$\\ \hline
5&(1,1,4) & $-13584186$&$2038284$\\ \hline
5&(3,0,4) & $-97309836 $&$13288212$\\ \hline
5&(0,3,3) & $-43113789$&$3446712$\\ \hline
5&(2,2,3) & $-278158212$&$30018033$\\ \hline
5&(4,1,3) & $-2127258450$&$ 237683160$\\ \hline
5&(6,0,3) & $-17629551342$&$1984125132$\\ \hline
5&(1,4,2) & $-\frac{1505585475}{2}$&$59337684$\\ \hline
5&(3,3,2) & $-5719059675 $&$538946784$\\ \hline
5&(5,2,2) & $-47352243843$&$4765188690$\\ \hline
5&(7,1,2) & $-416074383666 $&$43330628592$\\ \hline
5&(9,0,2) & $-3826973405736$&$407674059804$\\ \hline
5&(0,6,1) & $-1749516894$&$118262025$\\ \hline
5&(2,5,1) & $-\frac{58224030795}{4}$&$1152035055$\\ \hline
5&(4,4,1) & $-123411045312$&$11006956674$\\ \hline
5&(6,3,1) & $-1101441692466$&$104895765531$\\ \hline
5&(8,2,1) & $-10285610571117$&$1020477653388$\\ \hline
5&(10,1,1) & $-99935930525508 $&$ 10198147299480$\\ \hline
5&(12,0,1) & $-1006928375152050$&$104869613909880$\\ \hline
5&(1,7,0) & $-\frac{269759628747}{8}$&$2441617101$\\ \hline
5&(3,6,0) & $-306334343151 $&$24628038822$\\ \hline
5&(5,5,0) & $-2827530502056$&$247068849834$\\ \hline
5&(7,4,0) & $-\frac{54059816716407}{2}$&$2496139936101$\\ \hline
5&(9,3,0) & $-267911657815050$&$25688164768974$\\ \hline
5&(11,2,0) & $-2749583878851825$&$270571180960692$\\ \hline
5&(13,1,0) & $-29185732286251029$&$2922525674615178$\\ \hline
5&(15,0,0) & $-320196437236815048$&$32395182807076992$\\ \hline
\end{tabular}
\end{table}
\newpage
\begin{table}[H]
\centering
\caption{$g=0,M_6^4$} \label{g=0,N=6,k=4}
\begin{tabular}{|c|l|l|l|l|l|}
\hline
d&(a,b,c)& $N^0_{d,a,b,c}$  \\
\hline
1&(0,0,1) & $96$\\ \hline
1&(1,1,0) & $416$\\ \hline
1&(3,0,0) & $736$\\ \hline
2&(0,1,1) & $5568$\\ \hline
2&(2,0,1) & $21120$\\ \hline
2&(1,2,0)&$65536$ \\ \hline
2&(3,1,0) & $249856$\\ \hline
2&(5,0,0) & $1050368 $\\ \hline
3&(1,0,2) & $534528$\\ \hline
3&(0,2,1) & $1572864$\\ \hline
3&(2,1,1) & $9142272$\\ \hline
3&(4,0,1) & $53127168$\\ \hline
3&(1,3,0) & $29884416$\\ \hline
3&(3,2,0) & $173998080$\\ \hline
3&(5,1,0) & $1098065920$\\ \hline
3&(7,0,0) & $7453751296$\\ \hline
4&(0,0,3) & $13160448$\\ \hline
4&(1,1,2) & $294912000$\\ \hline
4&(3,0,2) & $2247081984$\\ \hline
4&(0,3,1) & $918552576$\\ \hline
4&(2,2,1) & $7183269888$\\ \hline
4&(4,1,1) & $56891523072$\\ \hline
4&(6,0,1) & $471610097664$\\ \hline
4&(1,4,0) & $24075304960$\\ \hline
4&(3,3,0) & $188432777216$\\ \hline
4&(5,2,0) & $1580113821696$\\ \hline
4&(7,1,0) & $14045898735616$\\ \hline
4&(9,0,0) & $132127907905536$\\ \hline
5&(0,1,3) & $8889827328$\\ \hline
5&(2,0,3) & $84506443776$\\ \hline
5&(1,2,2) & $263901413376$\\ \hline
5&(3,1,2) & $2562163605504$\\ \hline
5&(5,0,2) & $25399993171968$\\ \hline
5&(0,4,1) & $848390258688$\\ \hline
5&(2,3,1) & $8334908325888$\\ \hline
5&(4,2,1) & $83350340763648$\\ \hline
5&(6,1,1) & $867295299895296$\\ \hline
5&(8,0,1) & $9412512794935296$\\ \hline
5&(1,5,0) & $28347455242240$\\ \hline
5&(3,4,0) & $278636624084992$\\ \hline
5&(5,3,0) & $2909478586155008$\\ \hline
5&(7,2,0) & $31928942775238656$\\ \hline
5&(9,1,0) & $367130295597662208$\\ \hline
5&(11,0,0) & $4420079855887319040$\\ \hline
\end{tabular}
\end{table}
\begin{table}[H]
\centering
\caption{$g=1,M_6^4$} \label{g=1,N=6,k=4}
\begin{tabular}{|c|l|l|l|l|l|}
\hline
d&(a,b,c)& $N^1_{d,a,b,c}$&$E_{d,a,b,c}$  \\
\hline
1&(0,1,0) & $-\frac{260}{3}$&$0$\\ \hline
1&(2,0,0) & $-\frac{460}{3}$&$0$\\ \hline
2&(1,0,1)&$-4864$ &$0$\\ \hline
2&(0,2,0) & $ -\frac{43744 }{3}$&$0$\\ \hline
2&(2,1,0) & $-\frac{167008 }{3} $&$0$\\ \hline
2&(4,0,0) & $-\frac{698720 }{3} $&$0$\\ \hline
3&(0,0,2) & $-132864$&$768$\\ \hline
3&(1,1,1) & $-2246144 $&$13440$\\ \hline
3&(3,0,1) & $-13061376$&$73344$\\ \hline
3&(0,3,0) & $-\frac{21873152}{3}$&$27136$\\ \hline
3&(2,2,0) & $-42461696$&$161280$\\ \hline
3&(4,1,0) & $-\frac{804326144 }{3}$&$960512$\\ \hline
3&(6,0,0) & $-\frac{5455746304 }{3}$&$5898240$\\ \hline
4&(0,1,2) & $-74569728 $&$2571264$\\ \hline
4&(2,0,2) & $-571607040 $&$18837504$\\ \hline
4&(1,2,1) & $-1818378240$&$55492608$\\ \hline
4&(3,1,1) & $-14434416640$&$422731776$\\ \hline
4&(5,0,1) & $-119608152064$&$3456921600$\\ \hline
4&(0,4,0) & $-\frac{18345644032}{3}$&$146182144$\\ \hline
4&(2,3,0) & $-\frac{143574781952}{3}$&$1153077248$\\ \hline
4&(4,2,0) & $-\frac{1204880285696}{3}$&$9562345472$\\ \hline
4&(6,1,0) & $-\frac{10713421205504}{3}$&$83229306880$\\ \hline
4&(8,0,0) & $-33606505377792$&$754787155968$\\ \hline
5&(1,0,3) & $-21040939008 $&$2005917696$\\ \hline
5&(0,2,2) & $-65329790976$&$6142722048$\\ \hline
5&(2,1,2) & $-638763810816$&$57498402816$\\ \hline
5&(4,0,2) & $-6343662845952$&$561661526016$\\ \hline
5&(1,3,1) & $-2081253425152 $&$175519629312$\\ \hline
5&(3,2,1) & $-20872275427328$&$1716412514304$\\ \hline
5&(5,1,1) & $-217273901776896$&$17729712807936$\\ \hline
5&(7,0,1) & $-2356423499333632$&$191826913124352$\\ \hline
5&(0,5,0) & $-\frac{21453605306368}{3}$&$517300158464$\\ \hline
5&(2,4,0) & $-\frac{210823503020032}{3}$&$5103810117632$\\ \hline
5&(4,3,0) & $-\frac{2201894952304640}{3}$&$53278992793600$\\ \hline
5&(6,2,0) & $-\frac{24166727101579264}{3}$&$582424512397312$\\ \hline
5&(8,1,0) & $-\frac{277963197651501056}{3} $&$6634027667816448$\\ \hline
5&(10,0,0) & $-1116531188802797568 $&$78493526034874368$\\ \hline
\end{tabular}
\end{table}

\newpage
\section*{Acknowledgment}
We would like to thank Prof. H. Suzuki for valuable discussions.


\begin{thebibliography}{99}
\bibitem{AM} P. S. Aspinwall, D. R. Morrison.
\newblock{\em Topological field theory and rational curves.}
\newblock{Comm. Math. Phys. 151 (1993), no. 2, 245--262.}
\bibitem{BCOV} M. Bershadsky,  S. Cecotti, H. Ooguri,; C. Vafa. \newblock{\em Holomorphic anomalies in topological field theories.} 
\newblock{Nuclear Phys. B405 (1993), no.2--3, 279--304.}
\bibitem{bott-tu} R. Bott, L.-W. Tu. \newblock{\em Differential Forms in Algebraic Topology.}
\newblock{Graduate Texts in Mathematics, vol. 82, Springer (1982). }
\newblock{ISBN: 978-0-387-90613-3, https://doi.org/10.1007/978-1-4757-3951-0} 
\bibitem{virasoro1} T. Eguchi, K. Hori, and Chuan-Sheng Xiong.
\newblock{\em Quantum cohomology and Virasoro algebra}
\newblock{Phys.Lett. B402 (1997), 71-80.}
{http://dx.doi.org/10.1016/S0370-2693(97)00401-2}
\bibitem{virasoro2} T. Eguchi, M. Jinzenji,  and Chuan-Sheng Xiong.
\newblock{\em Quantum Cohomology and Free Field Representation}
\newblock{Nucl.Phys. B510 (1998), 608-622}
\bibitem{Getzler} E. Getzler.
\newblock{\em Intersection theory on $\overline {\mathcal{M}}_{1,4}$ and elliptic Gromov-Witten invariants}
\newblock{J. Amer. Math. Soc. 10 (1997), 973-998}
{https://doi.org/10.1090/S0894-0347-97-00246-4}
\bibitem{JK}M. Jinzenji, K. Kuwata.
\newblock{\em Elliptic Virtual Structure Constants and Generalizations of BCOV-Zinger Formula to Projective Fano
Hypersurfaces }
\newblock{J. High Energ. Phys. 2024, 135 (2024).}
\newblock{https://doi.org/10.1007/JHEP12(2024)135 }
\bibitem{JNS} M. Jinzenji, I. Nakamura, Y. Suzuki. {\em Conics on a Generic Hypersurface.}
\newblock{Preprint, arXiv:math/0412527. }
\bibitem{JR} M. Jinzenji.  
\newblock{\em Text copy of Maple worksheet for enumeration of degree $3$ elliptic curves in quintic hypersurface in $CP^5$ }
\newblock{.   Masao Jinzenji's ResearchGate homepage. https://www.researchgate.net/profile/Masao-Jinzenji\\   DOI: 10.13140/RG.2.2.32868.23680}
\bibitem{JS} M. Jinzenji, M. Shimizu
\newblock{\em Multi-point virtual structure constants and mirror computation of $CP^2$-model.}
\newblock{Commun. Number Theory Phys. 7 (2013), no. 3, 411--468. }
\newblock{https://dx.doi.org/10.4310/CNTP.2013.v7.n3.a2 }
\bibitem{Katz} S. Katz. {\em  Rational curves on Calabi-Yau manifolds: verifying predictions of Mirror Symmetry.}
\newblock{Preprint, arXiv:alg-geom/9301006.} 
\bibitem{Kont} M. Kontsevich. 
\newblock{\em Enumeration of rational curves via torus actions.}
\newblock{The moduli space of curves (Texel Island, 1994), 335--368.
Progr. Math., 129, Birkhäuser Boston, Inc., Boston, MA, 1995. }
\newblock{ISBN:0-8176-3784-2}
\bibitem{Kuwata}K. Kuwata. \newblock{\em Text copies of Mathematica files for "Enumeration of Elliptic Curves via Elliptic Gromov-Witten Invariants of Four Dimensional Projective Fano Hypersurfaces". }\\
\newblock{Ken Kuwata's ResearchGate homepage. https://www.researchgate.net/profile/Ken-Kuwata\\
DOI: 10.13140/RG.2.2.26180.28801}
\bibitem{KM}M. Kontsevich, Yu. Manin.\\
\newblock{\em Gromov-Witten classes, quantum cohomology,
and enumerative geometry.}
\newblock{Commun.Math. Phys. 164, (1994).}
\newblock{https://doi.org/10.1007/BF02101490}
\bibitem{KP} A. Klemm, R. Pandharipande.
\newblock{\em Enumerative geometry of Calabi-Yau 4-folds}
\newblock{\em Commun. Math. Phys. 281, (2008).}
\newblock{https://doi.org/10.1007/s00220-008-0490-9}
\bibitem{P1} R. Pandharipande.
\newblock{\em Hodge Integrals and Degenerate Contributions.}
\newblock{\em Comm Math Phys 208,  (1999).}
{https://doi.org/10.1007/s002200050766}
\end{thebibliography}
\end{document}